\documentclass[psreqno,reqno,11pt]{amsart}
\usepackage{eurosym}
\usepackage{amssymb}
\usepackage{tikz}
\usepackage{amsmath}
\usepackage{tikz}
\usepackage{pgflibraryplotmarks}
\usepackage{wasysym}
\usepackage{stmaryrd}

\theoremstyle{plain}

\theoremstyle{definition}

\theoremstyle{remark}

\begin{document}
\def\sect#1{\section*{\leftline{\large\bf #1}}}
\def\th#1{\noindent{\bf #1}\bgroup\it}
\def\endth{\egroup\par}

\title[Complex Vector Lattices]{%
Complex Vector Lattices Via Functional Completions}
\author{G. Buskes}
\address{Department of Mathematics, University of Mississippi, University,
MS 38677, USA}
\email{mmbuskes@olemiss.edu}
\author{C. Schwanke}
\address{Department of Mathematics, University of Mississippi, University,
MS 38677, USA}
\email{cmschwan@olemiss.edu}
\date{\today}

\begin{abstract}
We show that the Fremlin tensor product $C(X)\bar{\otimes}C(Y)$ is not square mean complete when X and Y are uncountable metrizable compact spaces. This motivates the definition of complexification of Archimedean vector lattices, the Fremlin tensor product of Archimedean complex vector lattices, and a theory of powers of Archimedean complex vector lattices.
\end{abstract}

\maketitle

\vspace{.5in}
\noindent
\textbf{1. Introduction} \vspace{.2in}

The standard references for the theory of vector lattices and Banach lattices
(see \cite{LuxZan1}, \cite{MN}, \cite{Schae3}, and \cite{Zan2}) all devote
some attention to complex vector lattices and complex Banach lattices, but a
reading of the treatment makes one feel that something is amiss. The
existence of a real cone in a complex vector space did show early promise,
and in fact, is essential at times in topics ranging from spectral theory
and vector measures to harmonic analysis. The emerging idea of a complex
modulus in the vector space complexification $E+iE$ of a Banach lattice $E$
dates to a 1963 paper by Rieffel (see \cite{Rieffel1} and also \cite{Rieffel2}) dealing with complex AL-spaces. In 1968 (see \cite{Lotz}), Lotz
defined, more generally, for Banach lattices $E$ and all $f,g\in E$ the
modulus $\left\vert f+ig\right\vert $ of an element $f+ig\in E+iE$ by%
\begin{equation*}
\left\vert f+ig\right\vert =\sup \{f\cos \theta +g\sin \theta :0\leq \theta
\leq 2\pi \}.\ (\ast )
\end{equation*}%
Luxemburg and Zaanen extend formula $(\ast )$ above to all uniformly
complete vector lattices (in \cite{LuxZan2}) in 1971, while studying order
bounded maps and integral operators. They realized that a theory of vector
lattices over $\mathbb{C}$ has to include a complex version of the
Kantorovich formula for the modulus of operators in the space of order
bounded operators $E\rightarrow F$, denoted by $\mathcal{L}_{b}(E,F)$, when $E$ is an Archimedean vector lattice and $F$ is Dedekind complete. That very
building block was provided in 1973 by de Schipper in \cite{dS}, with the
existence of the supremum in $(\ast )$ as a condition on $E$ and
Dedekind completeness of $F$ as follows. By defining a space of
complex order bounded operators $\mathcal{L}_{b}(E+iE,F+iF)$, de Schipper
proved that 
\begin{equation*}
\mathcal{L}_{b}(E,F)+i\mathcal{L}_{b}(E,F)=\mathcal{L}_{b}(E+iE,F+iF).
\end{equation*}%
Using the subscript $\mathbb{C}$ for the complexification of a vector space,
he thus proved that 
\begin{equation*}
\mathcal{L}_{b}(E_{\mathbb{C}},F_{\mathbb{C}})=\mathcal{L}_{b}(E,F)_{\mathbb{C}}\text{.}
\end{equation*}%
Interestingly, Luxemburg and Zaanen had proved the complex Kantorovich formula in the
earlier paper \cite{LuxZan2}, mentioned above, under the stronger condition
that $E$ is uniformly complete. Schaefer, in his book (\cite{Schae3}), defines
complex vector lattices axiomatically and derives formula $(\ast )$, but
includes uniform completeness in the axioms as well.

In spite of the validity of de Schipper's theorem under the mere
assumption of $(\ast )$, the assumption of uniform completeness
has proliferated in studies on complex vector lattices, almost invariably
identified with complexifications $E+iE$ of uniformly complete vector
lattices $E$. The choice of definition for complex vector lattices in \cite{Schae3} as well as the standard assumption of uniform completeness in
results for complex vector lattices in \cite{Zan2} appears to have codified
that practice. 

However, an alternative does exist in the literature, though it has hardly been used. Indeed, Mittelmeyer and Wolff in 1974 (see \cite{MW}) define what we call Archimedean complex vector lattices by axiomatizing an Archimedean modulus and they show
that the resulting Archimedean complex vector lattices are exactly the ones
that are vector space complexifications of Archimedean vector lattices with
property $(\ast )$. In light of the history sketched above, their complex
Archimedean vector lattices provide a ready made utility. The reader might
well ask: Why then write this paper?

One answer simply is this. Rewriting all of the theory for results that are
valid in Archimedean real vector lattices and Archimedean complex vector
lattices alike, seems a rather Herculean and, at times, uninteresting task.
We hasten to add that fundamental results for real vector lattices exist
that are not valid for complex ones. An example is the Riesz decomposition
property (see \cite{Vuz}). In the opposite direction, Kalton recently (see 
\cite{Kal}) proved surprising results for complex Banach
lattices that fail for real Banach lattices. In between there is a large
body of results that both theories have in common. But, even with complex
vector lattices satisfactorily defined in \cite{MW}, these results that are
in common, lack a proper transfer mechanism, a more or less
mechanical procedure that transfers real results into their complex
analogues, like de Schipper's result above.

In this paper, we present exactly such a mechanism. We do this in three ways.
First, we construct a vector lattice complexification for every Archimedean real
vector lattice, moving away from the vector space complexification for which
one needs to know a priori that one deals with a vector lattice in which
formula $(\ast )$ is valid. Secondly, these new complexifications are
precisely the Archimedean complex vector lattices introduced by Mittelmeyer
and Wolff. Thirdly, we show that these newly constructed complexifications
satisfy a natural universal property which, in many instances, tremendously
facilitates the transfer mechanism from real results to complex results. We introduce this vector lattice
complexification with a purpose in mind: differentiation
in Archimedean complex vector lattices via multilinear maps and tensor products. The
real version of such differentiation in vector lattices was introduced by
Loane in \cite{Loane}. A rapid development of polynomials on vector lattices
is currently under way and complex tensor products and complex powers of
complex vector lattices are needed. Motivated initially by this attempt to
complex differentiation, we started by looking at the real Fremlin tensor
product $E\bar{\otimes}E$ and were willing to assume uniform completeness of 
$E$, which has been the modus operandi in the literature, in order for 
\begin{equation*}
E\bar{\otimes}E+i(E\bar{\otimes}E)
\end{equation*}%
to be a complex vector lattice. However, $C(X)\bar{\otimes}C(Y)$ for any
uncountable compact metrizable spaces $X$ and $Y$ fails to have property $(\ast )
$ above. Introspection of the meaning of this failure naturally leads to the
concept of a new vector lattice complexification to address this deficiency.

This brings us to the content of this paper.

We indeed adopt the notion of a modulus on a vector space over $\mathbb{R}$
or $\mathbb{C}$, as introduced by Mittelmeyer and Wolff, in order to have
uniformity in language for results that are valid for both real and complex
vector lattices.

We have focused in this paper on results that are valid for
complex and real vector lattices. Indeed, we use the square mean
completion and its close ally, the vector lattice complexification to which
we alluded above, to construct a variety of Archimedean vector lattices
over $\mathbb{R}$ or $\mathbb{C}$, including the Fremlin tensor
product of Archimedean complex vector lattices, powers of Archimedean complex vector lattices, and Archimedean complex
vector lattices of maps of order bounded variation. These spaces in turn generalize a
host of results known for vector lattices over $\mathbb{R}$ to vector
lattices over $\mathbb{R}$ or $\mathbb{C}$, including a generalization of de
Schipper's result to multilinear maps of order bounded variation.

Finally, from our introduction above, the reader can rightly infer that a
literature search will find instances where uniform completeness was used by
habit rather than necessity to employ the vector space complexification in results on complex
vector lattices. To find such instances has not been a purpose of the
present paper. Instead, the methods in this paper provide an avenue to begin
a more systematic transfer from the vast literature on Archimedean real
vector lattices to the, in comparison, meagre set of results for Archimedean
complex vector lattices.

The first author gladly acknowledges a conversation with Arnoud van Rooij,
some twenty years ago, in which the idea and potential for a vector lattice
complexification were first raised.

\vspace{.5in}
\noindent
\textbf{2. Preliminaries} \vspace{.2in}

For all unexplained terminology about vector lattices we refer the reader to
the standard texts \cite{AB}, \cite{LuxZan1}, and \cite{Zan2}. Throughout, $\mathbb{K}$ stands for either $\mathbb{R}$ or $\mathbb{C}$, whereas $\mathbb{N}$ stands for the (nonzero) positive integers. For $s\in\mathbb{N}
$, we write $\times _{k=1}^{s}A_{k}$ for the Cartesian product $A_{1}\times
\dots \times A_{s}$, while $A\times \dots \times A\ $($s\ $times) is denoted
by $\times _{s}A$.

For the definition of an Archimedean vector lattice over $\mathbb{K}$,
central to this paper, we need the notion of a modulus on a vector space
(Mittelmeyer and Wolff in \cite{MW}).

\smallskip
\noindent
\th{Definition 2.1.} A modulus on
a vector space $E$ over $\mathbb{K}$ is an idempotent
mapping $m$ on $E$ that satisfies

\begin{enumerate}
\item[(1)] $m(\alpha f)=|\alpha |m(f)$ for every $\alpha
\in \mathbb{K}$ and for every $f\in E$,

\item[(2)] $m\Bigl(m\bigl(m(f)+m(g)\bigr)-m(f+g)\Bigr)
=m(f)+m(g)-m(f+g) $ for every $f,g\in E$, and

\item[(3)] $E$ is in the $\mathbb{K}$-linear hull
of $m(E) $.
\end{enumerate}
\endth

\noindent A modulus $m$ is said to be Archimedean
if for $f,g\in E$ it follows from $m\bigl(m(g)-nm(f)\bigr)=m(g)-nm(f)$ for every $n\in\mathbb{N}$ that $f=0$.

\smallskip We summarize the facts, obtained by Mittelmeyer and Wolff,
in the following theorem.

\smallskip
\noindent
\th{Theorem 2.2.}\ \textnormal{(\cite{MW}, Lemma 1.2, Corollary 1.4, Proposition 1.5, Theorem 2.2)}

\begin{enumerate}
\item[(1)] If $m$ is a modulus on a vector space $E$ over $\mathbb{R}$ then $m(E)$ is a cone in $E$. Moreover, $E$ is a vector lattice \textnormal{(}as defined in \cite{AB} and \cite{LuxZan1}\textnormal{)} under the partial ordering induced by $m(E)$. Furthermore, $m(f)=f\vee(-f)$ for every $f\in E$.

\item[(2)] If $m$ is an Archimedean modulus on a
vector space $E$ over $\mathbb{C}$ then $E$ is
of the form $E_{\rho }\oplus iE_{\rho }$, where $E_{\rho}:=m(E)-m(E)$ is an Archimedean vector lattice under the partial
ordering induced by $m(E)$. Moreover, $\sup \{f\cos \theta +g\sin
\theta :\theta \in \lbrack 0,2\pi ]\}$ exists in $E_{\rho }$ for every $f,g\in E_{\rho }$. Also, $m(f+ig)=\sup\{f\cos \theta +g\sin \theta :\theta \in[ 0,2\pi]\}$ for
every $f+ig\in E$.
\end{enumerate}
\endth

We call a vector space over $\mathbb{K}$ that is equipped with an Archimedean
modulus $m$ an \textit{Archimedean vector lattice over $\mathbb{K}$}. Since at times we will prove results for Archimedean vector lattices
over $\mathbb{K}$ that are known for Archimedean vector lattices over $%
\mathbb{R}$, we will use absolute value signs rather than $m$, for
convenience. The latter is justified by Theorem 2.2 above and, in turn, it
enables us from here on to talk about an Archimedean vector lattice
over $\mathbb{K}$ while suppressing the modulus $m$.

Let $E$ be an Archimedean vector lattice over $\mathbb{K}$. We
define $E^{+}:=\{f\in E:|f|=f\}$ and call $E^{+}$ the \textit{positive cone
of} $E$. Then $E_{\rho}=E^{+}-E^{+}$ by Theorem 2.2. A subset $A$ of
an Archimedean vector lattice $E$ over $\mathbb{K}$ is said to be \textit{%
order bounded} if there exists $f\in E^{+}$ such that $|a|\leq f$ for every $%
a\in A$. If $\sup A$ exists in $E^{+}$ for every nonempty order
bounded subset $A$ of $E^{+}$, we call $E$ \textit{Dedekind complete}. If $E$
is an Archimedean vector lattice over $\mathbb{R}$, we denote by $E^{\delta}$ the Dedekind completion of $E$.

Given Archimedean vector lattices $E$ and $F$ over $\mathbb{K}$, a $\mathbb{K}$-linear map $T:E\rightarrow F$ is called a \textit{vector lattice
homomorphism} if $T(|f|)=|T(f)|$ for every $f\in E$. A bijective vector
lattice homomorphism is called a \textit{vector lattice isomorphism}. If
there exists a vector lattice isomorphism between $E$ and $F$ we say that $E$
and $F$ are \textit{isomorphic as vector lattices over} $\mathbb{K}$.

More generally, let $V$ be a vector space over $\mathbb{R}$. We call 
$V$ the \textit{real part} of the complex vector space $V\oplus iV$. In the
latter case, we write $V=\left(V\oplus iV\right) _{\rho }$, in accordance
with the notation in Mittelmeyer and Wolff's Theorem 2.2 above. The complex vector space $V_{\mathbb{C}}:=V\oplus iV$ is called the \textit{vector space
complexification} of $V$. As usual, we consider $V$ to be a subset of $V_{\mathbb{C}}$ via the natural embedding. Given vector spaces $V_{1},...,V_{s},W$ over $\mathbb{R}$ and a map 
$T:\times _{k=1}^{s}V_{k\mathbb{C}}\rightarrow W_{\mathbb{C}}$, we say that $T$\textit{\ is real} if $T(f_{1},...,f_{s})\in W$ whenever $f_{k}\in V_{k}\
(k\in \{1,...,s\})$. The following formula (see Theorem 3 in \cite{BS}) uniquely extends
a map $T:\times_{k=1}^{s}V_{k}\rightarrow W$ which is linear over $\mathbb{R}$ in each variable separately (for short, $s_{\mathbb{R}}$\textit{-linear})
to a map $T_{\mathbb{C}}:\times_{k=1}^{s}V_{k\mathbb{C}}\rightarrow W_{\mathbb{C}}$ which is linear over $\mathbb{C}$ in each variable separately
(for short, $s_{\mathbb{C}}$\textit{-linear}):
\begin{equation*}
T_{\mathbb{C}}(f_{0}^{1}+if_{1}^{1},...,f_{0}^{s}+if_{1}^{s})=\sum\limits_{\epsilon_{k}\in \{0,1\}}T(f_{\epsilon _{1}}^{1},...,f_{\epsilon
_{s}}^{s})i^{\sum\limits_{k=1}^{s}\epsilon _{k}}
\end{equation*}
where $(f_{0}^{1}+if_{1}^{1},...,f_{0}^{s}+if_{1}^{s})\in \times
_{k=1}^{s}V_{k\mathbb{C}}$. We will say that $T_{\mathbb{C
}}$ is the \textit{complexification of }$T$. Conversely, when $T:\times
_{k=1}^{s}V_{k\mathbb{C}}\rightarrow W_{\mathbb{C}}$ is a, not necessarily
separately linear, real map then we write $T_{\rho}$ for the restriction of 
$T$ to $\times _{k=1}^{s}V_{k}$. It follows that $\left( T_{\rho }\right)_{\mathbb{C}}=T$ whenever $T$ is a real $s_{\mathbb{C}}$-linear map. We point
out (see \cite{Schae2}) that every vector space over $\mathbb{C}$ can be written as $V\oplus iV$ for some vector space $%
V $ over $\mathbb{R}$. An Archimedean vector lattice over $\mathbb{C}$,
however, contains a canonical real part that is determined by its modulus.
This fact has a variety of consequences because much of the basic
theory of Archimedean vector lattices over $\mathbb{R}$ is encoded via
vector lattice homomorphisms and positive linear maps and runs parallel with the theory
of Archimedean vector lattices over $\mathbb{C}$. We collect some examples
of this phenomenon that will be used repeatedly.

Let $E_{1},...,E_{s},F$ be Archimedean vector lattices over $%
\mathbb{K}$. A map $T:\times_{k=1}^{s}E_{k}\rightarrow F$ is called \textit{%
positive} if $T(f_{1},...,f_{s})\in F^{+}$ whenever $f_{k}\in E_{k}^{+}\ $%
for all $k\in \{1,...,s\}$. An $s_{\mathbb{K}}$\textit{-linear map} $%
T:\times _{k=1}^{s}E_{k}\rightarrow F$ is a map which is linear over $%
\mathbb{K}$ in each variable separately. An $s_{\mathbb{K}}$-linear map $T:\times _{k=1}^{s}E_{k}\rightarrow F$ is called a \textit{%
vector lattice} $s$-\textit{morphism} if for each $k\in \{1,...,s\}$ the
map from $E_{k}$ to $F$ defined by $f_{k}\mapsto T(f_{1},...,f_{k},...,f_{s})$ is a vector lattice homomorphism for fixed $f_{j}\in E_{j}^{+}\ (j\neq k)$. Every vector lattice $s$-morphism is
positive, and every positive $s_{\mathbb{K}}$-linear map is real. At times, for emphasis, we will refer to a vector lattice $s$-morphism between Archimedean vector lattices over $\mathbb{C}$, (respectively, Archimedean vector lattices over $\mathbb{R}$) as a vector lattice $s_{\mathbb{C}}$-morphism, or a vector lattice $\mathbb{C}$-homomorphism when $s=1$ (respectively, a vector lattice $s_{\mathbb{R}}$-morphism, or a vector lattice $\mathbb{R}$-homomorphism when $s=1$).

Given an Archimedean vector lattice $E$ over $\mathbb{K}$, we call a
vector subspace $L$ of $E$ a \textit{vector sublattice} of $E$ if $|f|\in L$ for every $f\in L$. If $L$ is a vector subspace of $E$ and for every $0<g\in E^{+}$ there exists $f\in L\cap E^{+}$ such that $0<f\leq g$, we say that $L$ is \textit{order dense} in $E$. If a vector subspace $L$ of $E$ has the property that $f\in L,g\in E$ and $|g|\leq|f|$ imply that $g\in L$, we call $L$ an \textit{ideal} in $E$. Hence every ideal of $E$ is a vector sublattice of $E$.

In \cite{BusSch}, we develop the uniform completion for Archimedean vector lattices over $\mathbb{R}$ in a manner that works just as well for developing the uniform completion for Archimedean vector lattices over $\mathbb{C}$. The definitions for relatively uniformly convergent sequences and relatively uniformly Cauchy sequences in an Archimedean vector lattice over $\mathbb{K}$, as well as the definition of a uniform completion of an Archimedean vector lattice over $\mathbb{K}$, are identical to what is found in Definitions 2.1 and 2.2 in \cite{BusSch}, modulo replacing $\mathbb{R}$ with $\mathbb{K}$. We also note that Propositions 3.1 and 3.2 in \cite{BusSch} also hold for Archimedean vector lattices over $\mathbb{K}$, and the proofs are identical to the real case. In particular, there exists an essentially unique uniform completion of every Archimedean vector lattice over $\mathbb{K}$.

\smallskip
\noindent
\th{Example 2.3.} If $E$ is a uniformly
complete Archimedean vector lattice over $\mathbb{R}$ then $E_{\mathbb{C}}$ is a
uniformly complete Archimedean vector lattice over $\mathbb{C}$. If $E$ is a uniformly complete Archimedean vector lattice over $\mathbb{C}$ then $E_{\rho}$ is a uniformly complete Archimedean vector lattice over $\mathbb{R}$. Dedekind complete Archimedean vector lattices over $\mathbb{K}$ are uniformly complete.
\endth

\smallskip
Let $E$ be an Archimedean vector lattice over $\mathbb{K}$ and let $A$ be a
subset of $E$. The \textit{pseudo uniform closure} $\bar{A}$ of $A$ is the
set of all $f\in E$ for which there exists a sequence $(f_{n})$ in $A$ such
that $f_{n}\overset{ru}{\rightarrow }f$ (also see page 85 of \cite{LuxZan1}). We call $A$ \textit{relatively uniformly closed} if $\bar{A}=A$, and we say that $A$ is \textit{uniformly dense} in $E$ if $\bar{A}=E$.

Like in \cite{BusSch}, we use transfinite induction to iterate the pseudo uniform closure of a vector sublattice $L$ of an Archimedean vector lattice $E$ over $\mathbb{K}$. For an Archimedean vector lattice $E$ over $\mathbb{K}$ and a nonempty subset $A$ of $E$, we define
\begin{equation*}
\begin{array}{l}
A_{1}:=A\text{,} \\ 
A_{\alpha}:=\overline{A_{\alpha-1}}\ \text{when}\ \alpha>1\ \text{is not a limit ordinal, and}\\ 
A_{\alpha}:=\bigcup_{\beta<\alpha}A_{\beta}\ \text{when}\ \alpha\ \text{is a limit ordinal.}
\end{array}
\end{equation*}

Given an Archimedean vector lattice $E$ over $\mathbb{K}$ and a nonempty subset $A\subseteq E$, we know from the (identical) complex analogue of Proposition 3.1 in \cite{BusSch} that $A_{\omega_{1}}$ is uniformly closed in $E$. It is quite possible however that there exists an ordinal $\alpha<\omega_{1}$ such that $A_{\alpha}$ is relatively uniformly closed in $E$. Motivated by this observation, we define the \textit{density number $\tau(A,E)$ of $A$ in $E$} by $\tau(A,E):=\min\lbrace\alpha:\bar{A}_{\alpha}=E\rbrace$.

\vspace{.5in}
\noindent
\textbf{3. Vector Lattice Complexifications} \vspace{.2in}

\smallskip
We discuss the specific case of Proposition 3.17 and Theorem 3.18 in \cite{BusSch} that we will use to complexify Archimedean vector lattices over $\mathbb{R}$ and multilinear maps over $\mathbb{R}$. Using the notation from Example 3.4 in \cite{BusSch}, let $\mu_{2,4}(x,y)=\sqrt{\frac{|x|^{2}+|y|^{2}}{2}}\ (x,y\in\mathbb{R})$. By Corollary 3.9 in \cite{BusSch} and the content following Corollary 3.10 in \cite{BusSch}, if $E$ is an Archimedean vector lattice over $\mathbb{R}$ and $f,g\in E$ then $\mu_{2,4}(f,g)=\frac{1}{\sqrt{2}}(f\boxplus g)$, where $f\boxplus g:=\sup\lbrace f\cos\theta+g\sin\theta:\theta\in[0,2\pi]\rbrace$, as defined by the authors of \cite{AzBoBus}. Therefore, from Mittelmeyer and Wolff's Theorem 2.2, a vector space $E+iE$ over $\mathbb{C}$ is an Archimedean vector lattice over $\mathbb{C}$ if and only if $E$ is a $\mu_{2,4}$-complete Archimedean vector lattice over $\mathbb{R}$. We refer to the $\mu_{2,4}$-completion $(E^{\mu_{2,4}},\phi)$ of $E$ as the \textit{square mean completion of $E$}. Noting that $\mu_{2,4}$ is absolutely invariant (defined following Theorem 3.15 in \cite{BusSch}), we summarize the newly found information regarding functional completions for this special case in the following corollary of Proposition 3.17 and Theorem 3.18 in \cite{BusSch}.

\smallskip
\noindent
\th{Corollary 3.1.} If $E$ is an Archimedean vector lattice over $\mathbb{R}$ then
there exists a unique square mean completion $(E^{\mu_{2,4}},\phi)$ of $E$. Moreover, if $E_{1},...,E_{s},F$ are
Archimedean vector lattices over $\mathbb{R}$ with square mean
completions $(E_{k}^{\mu_{2,4}},\phi_{k})\ (k\in \{1,...,s\})$ and $F$ is square mean complete, then for every vector lattice $s$-morphism $T:\times_{k=1}^{s}E_{k}\rightarrow F$ there exists a unique vector lattice $s$-morphism $T^{\mu_{2,4}}:\times
_{k=1}^{s}E_{k}^{\mu_{2,4}}\rightarrow F$ such that $T^{\mu_{2,4}}(\phi _{1}(f_{1}),...,\phi _{s}(f_{s}))=T(f_{1},...,f_{s})$. Furthermore, if $F$ is uniformly complete and $T:\times_{k=1}^{s}E_{k}\rightarrow F$ is a positive $s$-linear map then there exists a unique positive $s$-linear map $T^{\mu_{2,4}}:\times_{k=1}^{s}E_{k}^{\mu_{2,4}}\rightarrow F$ such that $T^{\mu_{2,4}}(\phi_{1}(f_{1}),...,\phi_{s}(f_{s})=T(f_{1},...,f_{s})$ for every $f_{k}\in E_{k}\ (k\in\lbrace 1,...,s\rbrace)$. Here $\phi_{k}$ is the natural embedding of $E_{k}$ into $E_{k}^{u}$.
\endth

\smallskip
We now turn to complexifications of Archimedean vector
lattices over $\mathbb{R}$.

\smallskip
\noindent
\th{Definition 3.2.} \textit{For an
Archimedean vector lattice }$E$\textit{\ over }$\mathbb{R}$\textit{\ we
define a pair }$(E_{|\mathbb{C}|},\phi )$\textit{\ to be a vector lattice
complexification of }$E$\textit{\ if the following hold.}

\begin{enumerate}
\item[(1)] $E_{|\mathbb{C}|}$\textit{\ is an Archimedean vector
lattice over }$\mathbb{C}$\textit{.}

\item[(2)] $\phi:E\rightarrow(E_{|\mathbb{C}|})_{\rho}$\textit{\ is an injective vector lattice $\mathbb{R}$-homomorphism.}

\item[(3)] \textit{For every Archimedean vector lattice} $F$\textit{over $\mathbb{C}$ as well as for every vector lattice $\mathbb{R}$-homomorphism} $T:E\rightarrow F_{\rho}$,\textit{\ there exists a unique vector lattice $\mathbb{C}$-homomorphism }$T_{|\mathbb{C}|}:E_{|\mathbb{C}|}\rightarrow F$\textit{\ such that }$T_{|\mathbb{C}|}\circ \phi =T$\textit{.}
\end{enumerate}
\endth

\smallskip
We next prove the existence and uniqueness of vector
lattice complexifications.

\smallskip
\noindent
\th{Theorem 3.3.} \textit{\ If }$E$\textit{\ is an Archimedean vector lattice over }$\mathbb{R} $\textit{\ then
there exists a vector lattice complexification of }$E$\textit{,
unique up to vector lattice isomorphism. }
\endth

\begin{proof}
Let $E$ be an Archimedean vector lattice over $\mathbb{R}$. By Corollary 3.1, there exists a unique square mean completion $(E^{\mu_{2,4}},\phi)$ of $E$. Define $E_{|\mathbb{C}|}:=(E^{\mu_{2,4}})_{\mathbb{C}}$ and observe that $E_{|\mathbb{C}|}$ is an Archimedean vector lattice over $\mathbb{C}$
and that $(E_{|\mathbb{C}|})_{\rho}=E^{\mu_{2,4}}$. Next, let $F$ be an Archimedean vector lattice over $\mathbb{C}$ and let $T:E\rightarrow F_{\rho}$
be a vector lattice $\mathbb{R}$-homomorphism. Since $F_{\rho}$ is square mean complete, there exists a unique vector lattice
$\mathbb{R}$-homomorphism $T^{\mu_{2,4} }:E^{\mu_{2,4}}\rightarrow F_{\rho}$ such that $T^{\mu_{2,4}}\circ\phi=T$. Define $T_{|\mathbb{C}|}:E_{|\mathbb{C}|}\rightarrow F$ by $T_{|\mathbb{C}|}(f+ig)=T^{\mu_{2,4}}(f)+iT^{\mu_{2,4}}(g)$ for every $f+ig\in E_{|\mathbb{C}|}$. Then $T_{|\mathbb{C}|}\circ\phi=T$. Moreover, for $f+ig\in E_{|\mathbb{C}|}$ we have from Corollary 3.13 in \cite{BusSch} (see also Proposition 3.4 of \cite{Az}) that
\begin{align*}
T_{|\mathbb{C}|}(|f+ig|)=T^{\mu_{2,4}}(f\boxplus g)=T^{\mu_{2,4}}(f)\boxplus T^{\mu_{2,4}}(g)=|T_{|\mathbb{C}|}(f+ig)|.
\end{align*}
\noindent
Thus $T_{|\mathbb{C}|}$ is a vector lattice $\mathbb{C}$-homomorphism and therefore $(E_{|\mathbb{C}|},\phi )$ is a vector lattice complexification of $E$. Next, we prove the uniqueness. To this end, suppose $(E_{1|\mathbb{C}|},\phi _{1})$ and $(E_{2|\mathbb{C}|},\phi _{2})$ are vector lattice complexifications of $E$. Then $((E_{1|\mathbb{C}|})_{\rho},\phi_{1})$ and $((E_{2|\mathbb{C}|})_{\rho },\phi_{2})$ are square mean completions of $E$, and hence there exists a vector lattice isomorphism $\gamma:(E_{1|\mathbb{C}|})_{\rho}\rightarrow(E_{2|\mathbb{C}|})_{\rho}$. Similar to $T_{|\mathbb{C}|}$ above, the map $\gamma_{\mathbb{C}}:E_{1|\mathbb{C}|}\rightarrow E_{2|\mathbb{C}|}$ defined by $\gamma_{\mathbb{C}}(f+ig)=\gamma(f)+i\gamma(g)$ is a vector lattice $\mathbb{C}$-homomorphism. The bijectivity of $\gamma_{\mathbb{C}}$ is evident.
\end{proof}

\smallskip
For the square mean completion $(E^{\mu_{2,4}},\phi)$ of $E$, we will from now on identify $E$ with $\phi(E)$. Using this identification, we complexify positive $s_{\mathbb{R}}$-linear maps (respectively, vector lattice $s_{\mathbb{R}}$-morphisms) to positive $s_{\mathbb{C}}$-linear maps (respectively, vector lattice $s_{\mathbb{C}}$-morphisms) as follows. Let $E_{1},...,E_{s},F$ be Archimedean vector lattices over $\mathbb{R}$ with $F$ square mean complete, and let $T:\times_{k=1}^{s}E_{k}\rightarrow F$ be a vector lattice $s_{\mathbb{R}}$-morphism. For $(f_{0}^{1}+if_{1}^{1},...,f_{0}^{s}+if_{1}^{s})\in\times_{k=1}^{s}E_{k|\mathbb{C}|}$, define $T_{|\mathbb{C}|}:\times_{k=1}^{s}E_{k|\mathbb{C}|}\rightarrow F_{\mathbb{C}}$ by
\begin{equation*}
T_{|\mathbb{C}|}(f_{0}^{1}+if_{1}^{1},...,f_{0}^{s}+if_{1}^{s}):=\sum\limits_{\epsilon_{k}\in\{0,1\}}T^{\mu_{2,4}}(f_{\epsilon _{1}}^{1},...,f_{\epsilon_{s}}^{s})i^{\sum\limits_{k=1}^{s}\epsilon_{k}}.
\end{equation*}

If $F$ is uniformly complete and $T$ above is any positive $s_{\mathbb{R}}$-linear map, we define $T_{|\mathbb{C}|}$ in a similar manner. We collect a few facts regarding this complexification in the following proposition. Statement (3) and the statement that $T_{|\mathbb{C}|}=\left(T^{\mu_{2,4}}\right)_{\mathbb{C}}$ in (1) and (2) are evident. The proof of (2) follows from Corollary 3.1, and the proof of (1) is similar to the complexification of vector lattice homomorphisms seen in the proof of Theorem 3.3.

\smallskip
\noindent
\th{Proposition 3.4.} Let $E_{1},...,E_{s},F$ be Archimedean vector lattices over $\mathbb{R}$ with $F$ square mean complete.
\begin{itemize}
\item[(1)] If a map $T:\times_{k=1}^{s}E_{k}\rightarrow F$ is a vector lattice $s_{\mathbb{R}}$-morphism then $T_{|\mathbb{C}|}$ is a vector lattice $s_{\mathbb{C}}$-morphism and $T_{|\mathbb{C}|}=\left(T^{\mu_{2,4}}\right)_{\mathbb{C}}$.

\item[(2)] If $F$ is uniformly complete and $T:\times_{k=1}^{s}E_{k}\rightarrow F$ is a positive $s_{\mathbb{R}}$-linear map then $T_{|\mathbb{C}|}$ is a positive $s_{\mathbb{C}}$-linear map and $T_{|\mathbb{C}|}=\left(T^{\mu_{2,4}}\right)_{\mathbb{C}}$.

\item[(3)] If in (1) or (2) all $E_{1},...,E_{s}$ are square mean complete then $T_{|\mathbb{C}|}=T_{\mathbb{C}}$.
\end{itemize}
\endth

\vspace{.5in}
\noindent
\textbf{4. The Archimedean Vector Lattice Tensor Product} \vspace{.2in}

In this section, we define the tensor product of
Archimedean vector lattices over $\mathbb{K}$ and prove the existence of the Archimedean complex tensor product by complexifying the Fremlin tensor product of Archimedean real vector lattices.

We start with the definition of these tensor products of
Archimedean vector lattices over $\mathbb{K}$. For $s=2$ and $\mathbb{K}=\mathbb{R}$ the
definition coincides with Fremlin's definition of the Archimedean tensor
product of Archimedean vector lattices over $\mathbb{R}$ (see \cite{Fr}).

\smallskip
\noindent
\th{Definition 4.1.} Given Archimedean vector
lattices $E_{1},...,E_{s}$ over $\mathbb{K}$, we define a pair $(\bar{%
\otimes }_{k=1}^{s}E_{k},\bar{\otimes })$ to be an Archimedean vector
lattice tensor product of $E_{1},...,E_{s}$ if the following hold.

\begin{enumerate}
\item[(1)] \textit{$\bar{\otimes }_{k=1}^{s}E_{k}$ is an Archimedean
vector lattice over $\mathbb{K}$. }

\item[(2)] \textit{$\bar{\otimes }$ is a vector lattice $s$-morphism.}

\item[(3)] \textit{For every Archimedean vector lattice $F$ over $\mathbb{K}$
and for every vector lattice $s$-morphism $T:\times
_{k=1}^{s}E_{k}\rightarrow F$, there exists a uniquely determined vector lattice
homomorphism $T^{\bar{\otimes }}:\bar{\otimes}_{k=1}^{s}E_{k}\rightarrow F$ such that $T^{\bar{\otimes }}\circ \bar{\otimes }=T$%
. }
\end{enumerate}
\endth

\smallskip
Following its original proof, one can extend Fremlin's Theorem 4.2 in \cite{Fr} to the vector lattice tensor product of any number of factors (see also \cite{{Sh}}, Section 2). Below and throughout the rest of this section, $(\otimes_{k=1}^{s}V_{k},\otimes)$ denotes the algebraic tensor product of vector spaces $V_{1},...,V_{s}$ over $\mathbb{K}$.

\smallskip
\th{Lemma 4.2.} Let $E_{1},...,E_{s}$ be Archimedean vector lattices over $\mathbb{R}$.

\begin{enumerate}
\item[(1)] \textit{There exists an essentially unique Archimedean vector lattice $\bar{\otimes}_{k=1}^{s}E_{k}$ over $\mathbb{R}$ \textit{\ and a vector lattice }$s$\textit{-morphism }$\bar{\otimes }:\times_{k=1}^{s}E_{k}\rightarrow\bar{\otimes}_{k=1}^{s}E_{k}$\textit{\ such
that for every Archimedean vector lattice }$F$ over $\mathbb{R}$ \textit{and every vector lattice }$s$\textit{-morphism }$T:\times
_{k=1}^{s}E_{k}\rightarrow F$\textit{, there exists a unique vector lattice
homomorphism }$T^{\bar{\otimes }}:\bar{\otimes }%
_{k=1}^{s}E_{k}\rightarrow F$\textit{\ such that }$T^{\bar{\otimes }%
}\circ \bar{\otimes }=T$\textit{. } }

\item[(2)] \textit{There exists an injective linear map $S:\otimes
_{k=1}^{s}E_{k}\rightarrow\bar{\otimes}_{k=1}^{s}E_{k}$
such that $S\circ\otimes=\bar{\otimes}$.}

\item[(3)] For every $w\in\bar{\otimes}_{k=1}^{s}E_{k}$\textit{, there exist }$x_{k}\in E_{k}^{+}\ (k\in \{1,...,s\})$\textit{\
such that for every }$\epsilon >0$\textit{, there exists }$v\in \otimes
_{k=1}^{s}E_{k}$\textit{\ such that }$|w-v|\leq \epsilon (x_{1}\otimes \dots
\otimes x_{s})$\textit{, i.e. }$\otimes _{k=1}^{s}E_{k}$\textit{\ is
relatively uniformly dense in }$\bar{\otimes }_{k=1}^{s}E_{k}$\textit{.}

\item[(4)] \textit{For every} $0<w\in\bar{\otimes}_{k=1}^{s}E_{k}$\textit{\ there exist }$x_{k}\in E_{k}^{+}\ (k\in
\{1,...,s\})$\textit{\ such that }$0<(x_{1}\otimes \dots \otimes x_{s})\leq
w$\textit{, i.e. }$\otimes _{k=1}^{s}E_{k}$\textit{\ is order dense in }$
\bar{\otimes}_{k=1}^{s}E_{k}$\textit{.}
\end{enumerate}
\endth

\smallskip
The main result of this section deals with the existence and uniqueness of the complex Archimedean vector lattice tensor product and requires several prerequisite results. The next lemma surely is known but we could only find an explicit reference in the literature for a special case in the thesis \cite{vZ}.

\smallskip
\noindent
\th{Lemma 4.3.} If $V_{1},...,V_{s}$ are vector
spaces over $\mathbb{R}$ then $\otimes _{k=1}^{s}(V_{k\mathbb{C}})$ and $(\otimes _{k=1}^{s}V_{k})_{\mathbb{C}}$ are isomorphic as vector spaces over 
$\mathbb{C}$.
\endth

\begin{proof} Since the algebraic tensor product is associative, we only need to prove the
result for $s=2$, and use induction. The case $s=2$ is the content of Theorem 2.1.2 in \cite{vZ}, but, we provide a sketch of van Zyl's proof to correct some potential confusion caused by an accumulation of minor misprints. First let $U$ and $V$
be vector spaces over $\mathbb{R}$, and let $(U\otimes V,\otimes )$ and $(U_{%
\mathbb{C}}\otimes _{1}V_{\mathbb{C}})$ be the algebraic tensor products of $%
U,V$, respectively $U_{\mathbb{C}},V_{\mathbb{C}}$. Since $\otimes _{\mathbb{%
C}}:U_{\mathbb{C}}\times V_{\mathbb{C}}\rightarrow (U\otimes V)_{\mathbb{C}}$
is a bilinear map over $\mathbb{C}$, it induces a unique $\mathbb{C}$-linear
map $T:U_{\mathbb{C}}\otimes _{1}V_{\mathbb{C}}\rightarrow (U\otimes V)_{%
\mathbb{C}}$. It is easy to see that $T$ is surjective. To show that $T$ is
injective, let $w=\sum\limits_{k=1}^{n}(u_{k}+iu_{k}^{\prime })\otimes
_{1}(v_{k}+iv_{k}^{\prime })\in U_{\mathbb{C}}\otimes _{1}V_{\mathbb{C}}$
and suppose that $T(w)=0$. Note that $T(w)=\sum\limits_{k=1}^{n}(u_{k}%
\otimes v_{k}-u_{k}^{\prime }\otimes v_{k}^{\prime }+iu_{k}^{\prime }\otimes
v_{k}+iu_{k}\otimes v_{k}^{\prime })$, and so for any $\mathbb{R}$-linear
functionals $\phi $ on $U$ and $\psi $ on $V$ we have
\begin{equation}
\sum\limits_{k=1}^{n}\bigl(\phi (u_{k})\psi (v_{k})-\phi (u_{k}^{\prime})\psi(v_{k}^{\prime })\bigr)=0\ \text{and}\ \sum\limits_{k=1}^{n}\bigl(\phi(u_{k}^{\prime})\psi(v_{k})+\phi(u_{k})\psi(v_{k}^{\prime})\bigr)=0\ \tag{$\ast$}.
\end{equation}

Let $\xi=\xi_{r}+i\xi_{c}$ be a $\mathbb{C}$-linear functional on $%
U_{\mathbb{C}}$ and let $\eta =\eta _{r}+i\eta _{c}$ be a $\mathbb{C}$%
-linear functional on $V_{\mathbb{C}}$, both written in their natural
decompositions. Then $\xi _{r},\xi _{c}$ are $\mathbb{R}$-linear functionals
on $U$ and $\eta _{r},\eta _{c}$ are $\mathbb{R}$-linear functionals on $V$.
Now

\begin{align*}
&\sum\limits_{k=1}^{n}\xi(u_{k}+iu^{\prime }_{k})\eta(v_{k}+iv^{\prime }_{k})
\\
=&\sum\limits_{k=1}^{n}\bigl(\xi_{r}(u)\eta_{r}(v_{k})-\xi_{r}(u^{\prime
}_{k})\eta_{r}(v^{\prime }_{k})\bigr)-\sum\limits_{k=1}^{n}\bigl(%
\xi_{r}(u^{\prime }_{k})\eta_{c}(v_{k})+\xi_{r}(u_{k})\eta_{c}(v^{\prime
}_{k})\bigr) \\
+&i\sum\limits_{k=1}^{n}\bigl(\xi_{r}(u^{\prime
}_{k})\eta_{r}(v_{k})+\xi_{r}(u_{k})\eta_{r}(v^{\prime }_{k})\bigr)%
+i\sum\limits_{k=1}^{n}\bigl(\xi_{r}(u_{k})\eta_{c}(v_{k})-\xi_{r}(u^{\prime
}_{k})\eta_{c}(v^{\prime }_{k})\bigr) \\
-&\sum\limits_{k=1}^{n}\bigl(\xi_{c}(u^{\prime
}_{k})\eta_{r}(v_{k})+\xi_{c}(u_{k})\eta_{r}(v^{\prime }_{k})\bigr)%
-\sum\limits_{k=1}^{n}\bigl(\xi_{c}(u_{k})\eta_{c}(v_{k})-\xi_{c}(u^{\prime
}_{k})\eta_{c}(v^{\prime }_{k})\bigr) \\
+&i\sum\limits_{k=1}^{n}\bigl(\xi_{c}(u_{k})\eta_{r}(v_{k})-\xi_{c}(u^{%
\prime }_{k})\eta_{r}(v^{\prime }_{k})\bigr)-i\sum\limits_{k=1}^{n}\bigl(%
\xi_{c}(u^{\prime }_{k})\eta_{c}(v_{k})+\xi_{c}(u_{k})\eta_{c}(v^{\prime
}_{k})\bigr).
\end{align*}

Applying $(\ast)$ again to each of these eight summands, we
have that $\sum\limits_{k=1}^{n}\xi (u_{k}+iu_{k}^{\prime })\eta
(v_{k}+iv_{k}^{\prime })=0$. Therefore $w=0$ and $T$ is injective. Then $T$
is a vector space isomorphism.
\end{proof}

\smallskip
In light of the previous lemma, we will from now identify $(\otimes_{k=1}^{s}V_{k\mathbb{C}})_{\rho}$ with $\otimes_{k=1}^{s}V_{k}$ for vector spaces $V_{1},...,V_{s}$ over $\mathbb{R}$.

Next, we note that there exists a simpler construction of the square mean completion than the construction preceding Proposition 3.16 in \cite{BusSch}, which was given in a more general setting. Indeed, in Remark 4 of \cite{Az}, Azouzi constructs a square mean completion of an Archimedean vector lattice $E$ over $\mathbb{R}$ essentially as follows. Let $E_{1}:=E$ and for every $n\in\mathbb{N}$, define $E_{n+1}:=E_{n}\cup[\lbrace\mu_{2,4}(f,g):f,g\in E_{n}\rbrace]$, where $[\lbrace\mu_{2,4}(f,g):f,g\in E_{n}\rbrace]$ denotes the vector subspace of $E^{\delta}$ generated by $\lbrace\mu_{2,4}(f,g):f,g\in E_{n}\rbrace$. Then define $E^{\boxplus}:=\bigcup_{n\in\mathbb{N}}E_{n}$. To see that $E^{\boxplus}$ is a vector lattice, note that for every $f\in E^{\boxplus}$ there exists $n\in\mathbb{N}$ such that $f\in E_{n}$. Then $|f|=\sqrt{2}\mu_{2,4}(f,0)\in E_{n+1}$. It follows that $E^{\boxplus}$ is the square mean completion of $E$, that is, $E^{\boxplus}$ and $E^{\mu_{2,4}}$ are isomorphic as vector lattices. In fact, from the identity $\lambda\mu_{2,4}(f,g)=\mu_{2,4}(\lambda f,\lambda g)$ for every $\lambda\in\mathbb{R}^{+}$ and every $f,g\in E$, we have $E_{n+1}^{+}=\lbrace\sum\limits_{k=1}^{m}\mu_{2,4}(f_{k},g_{k}):f_{k},g_{k}\in E_{n}\rbrace$. We use this fact in the first of the two following lemmas that are needed for Proposition 4.6.

\smallskip
\noindent
\th{Lemma 4.4.} Denote the standard sine and cosine functions on $[0,\frac{\pi}{2}]$ by $\sin$ and $\cos$, respectively. For an Archimedean vector lattice $E$ over $\mathbb{R}$ and for every $f\in(E^{\mu_{2,4}})^{+}$ there exists $u_{1},...,u_{n}\in E^{+}$ and $t_{k,1},...,t_{k,p_{k}}\in\lbrace\cos,\sin\rbrace\ (k\in\lbrace 1,...,n\rbrace)$ such that 
\begin{equation*}
f=\underset{\theta_{k}\in[0,\frac{\pi}{2}]}{\sup}\Bigl\lbrace\sum\limits_{k=1}^{n}\prod\limits_{j=1}^{p_{k}}t_{k,j}(\theta_{k})u_{k}\Bigr\rbrace.
\end{equation*}
\endth

\begin{proof}
Our proof is via mathematical induction. Let $h\in E_{2}^{+}$ and first suppose  that $f=\mu_{2,4}(u,v)$ for some $u,v\in E^{+}$. Then $f=\sup\lbrace u\cos\theta+v\sin\theta:\theta\in[0,\frac{\pi}{2}]\rbrace$. Next, suppose that $f=\sum\limits_{k=1}^{n}\mu_{2,4}(u_{k},v_{k})$. Then
\begin{equation*}
f=\sum\limits_{k=1}^{n}\underset{\theta_{k}\in[0,\frac{\pi}{2}]}{\sup}\lbrace u_{k}\cos\theta_{k}+v_{k}\sin\theta_{k}\rbrace=\underset{\theta_{k}\in[0,\frac{\pi}{2}]}{\sup}\Bigl\lbrace\sum\limits_{k=1}^{n}(u_{k}\cos\theta_{k}+v_{k}\sin\theta_{k})\Bigr\rbrace.
\end{equation*}
\noindent
This completes the base step of the induction argument. For the inductive step, suppose that for every $f\in E_{n}^{+}$ there exists $u_{1},...,u_{n}\in E^{+}$ and $t_{1},...,t_{p_{k}}\in\lbrace\cos,\sin\rbrace\ (k\in\lbrace 1,...,n\rbrace)$ such that
\begin{equation*}
f=\underset{\theta_{k,j}\in[0,\frac{\pi}{2}]}{\sup}\Bigl\lbrace\sum\limits_{k=1}^{n}\prod\limits_{j=1}^{p_{k}}t_{k,j}(\theta_{k,j})u_{k}\Bigr\rbrace.
\end{equation*}
\noindent
Let $f\in E_{n+1}^{+}$. From the argument in the base step above, we may assume that $f=\mu_{2,4}(u,v)$ for some $u,v\in E_{n}^{+}$. By the induction hypothesis there exists $u_{1},...,u_{n},v_{1},...,v_{n}\in E^{+}$ and $t_{k,1},...,t_{k,p_{k}},s_{k,1},...,s_{k,r_{k}}\in\lbrace\cos,\sin\rbrace\ (k\in\lbrace 1,...,s\rbrace)$ such that
\begin{equation*}
u=\underset{\theta_{k,j}\in[0,\frac{\pi}{2}]}{\sup}\Bigl\lbrace\sum\limits_{k=1}^{n}\prod\limits_{j=1}^{p_{k}}t_{k,j}(\theta_{k,j})u_{k}\Bigr\rbrace\ \text{and}\ v=\underset{\theta_{k,j}\in[0,\frac{\pi}{2}]}{\sup}\Bigl\lbrace\sum\limits_{k=1}^{m}\prod\limits_{j=1}^{r_{k}}s_{k,j}(\theta_{k,j})v_{k}\Bigr\rbrace.
\end{equation*}
\noindent
Then
\begin{align*}
f&=\mu_{2,4}\Bigl(\underset{\theta_{k,j}\in[0,\frac{\pi}{2}]}{\sup}\Bigl\lbrace\sum\limits_{k=1}^{n}\prod\limits_{j=1}^{p_{k}}t_{k,j}(\theta_{k,j})u_{k}\Bigr\rbrace,\underset{\theta_{k,j}\in[0,\frac{\pi}{2}]}{\sup}\Bigl\lbrace\sum\limits_{k=1}^{m}\prod\limits_{j=1}^{r_{k}}s_{k,j}(\theta_{k,j})v_{k}\Bigr\rbrace\Bigr)\\
&=\underset{\phi\in[0,\frac{\pi}{2}]}{\sup}\biggl\lbrace\underset{\theta_{k,j}\in[0,\frac{\pi}{2}]}{\sup}\Bigl\lbrace\sum\limits_{k=1}^{n}\prod\limits_{j=1}^{p_{k}}t_{k,j}(\theta_{k,j})u_{k}\Bigr\rbrace\cos\phi+\underset{\theta_{k,j}\in[0,\frac{\pi}{2}]}{\sup}\Bigl\lbrace\sum\limits_{k=1}^{m}\prod\limits_{j=1}^{r_{k}}s_{k,j}(\theta_{k,j})v_{k}\Bigr\rbrace\sin\phi\biggr\rbrace\\
&=\underset{\phi,\theta_{j,k}\in[0,\frac{\pi}{2}]}{\sup}\Bigl\lbrace\sum\limits_{k=1}^{n}\prod\limits_{j=1}^{p_{k}}t_{k,j}(\theta_{k,j})\cos\phi u_{k}+\sum\limits_{k=1}^{m}\prod\limits_{j=1}^{r_{k}}s_{k,j}(\theta_{k,j})\sin\phi v_{k}\Bigr\rbrace.
\end{align*}
\noindent
\end{proof}

\smallskip
The next lemma can be verified using mathematical induction. We do not include the proof.

\smallskip
\noindent
\th{Lemma 4.5.} Let $t_{1},...,t_{n}$ be Lipschitz functions on $\mathbb{R}$ with Lipschitz constant $1$. Also assume that $|t_{k}(x)|\leq 1$ for every $k\in\lbrace 1,...,n\rbrace$ and every $x\in\mathbb{R}$. Then for every $x_{k},y_{k}\in\mathbb{R}\ (k\in\lbrace 1,...,n\rbrace)$ we have
$\bigl|\prod\limits_{k=1}^{n}t_{k}(x_{k})-\prod\limits_{k=1}^{n}t_{k}(y_{k})\bigr|\leq\sum\limits_{k=1}^{n}|x_{k}-y_{k}|$.
\endth

\smallskip
The idea of the proof for the following proposition comes from Lemma 2.8 in \cite{Az}.

\smallskip
\noindent
\th{Proposition 4.6.} If $E$ is an Archimedean vector lattice over $\mathbb{R}$ then $E$ is relatively uniformly dense in $E^{\mu_{2,4}}$.
\endth

\begin{proof}
Let $E$ be an Archimedean vector lattice over $\mathbb{R}$ and first suppose that $f\in(E^{\mu_{2,4}})^{+}$. Say that $f=\underset{\theta_{k,j}\in[0,\frac{\pi}{2}]}{\sup}\lbrace\sum\limits_{k=1}^{n}\prod\limits_{j=1}^{p_{k}}t_{k,j}(\theta_{k,j})u_{k}\rbrace$ for some $u_{1},...,u_{n}\in E^{+}$ and $t_{k,1},...,t_{k,p_{k}}\in\lbrace\cos,\sin\rbrace\ (k\in\lbrace 1,...,n\rbrace)$. Note that given $\theta_{k,j}\in[0,\frac{\pi}{2}]$ $m\in\mathbb{N}$ there exist $l_{k,j}\in\mathbb{N}$ such that $|\frac{l_{k,j}\pi}{2^{m}}-\theta_{k,j}|\leq\frac{\pi}{2^{m}}$. Since sine and cosine are both Lipschitz functions with Lipschitz constant $1$ we have from Lemma 4.5 that
\begin{align*}
\Bigl|\sum\limits_{k=1}^{n}\prod\limits_{j=1}^{p_{k}}t_{k,j}(\theta_{k,j})u_{k}-\sum\limits_{k=1}^{n}\prod\limits_{j=1}^{p_{k}}t_{k,j}(\frac{l_{k,j}\pi}{2^{m}})u_{k}\Bigr|&\leq\sum\limits_{k=1}^{n}\Bigl|\prod\limits_{j=1}^{p_{k}}t_{k,j}(\theta_{k,j})-\prod\limits_{j=1}^{p_{k}}t_{k,j}(\frac{l_{k,j}\pi}{2^{m}})\Bigr||u_{k}|\\
&\leq\sum\limits_{k=1}^{n}\sum\limits_{j=1}^{p_{k}}|\theta_{k,j}-\frac{l_{k,j}\pi}{2^{m}}||u_{k}|\\
&\leq\frac{\pi}{2^{m}}\sum\limits_{k=1}^{n}p_{k}|u_{k}|.
\end{align*}

\smallskip
Thus,
\begin{align*}
\sum\limits_{k=1}^{n}\prod\limits_{j=1}^{p_{k}}t_{k,j}(\theta_{k,j})u_{k}&\leq\sum\limits_{k=1}^{n}\prod\limits_{j=1}^{p_{k}}t_{k,j}(\frac{l_{k,j}\pi}{2^{m}})u_{k}+\frac{\pi}{2^{m}}\sum\limits_{k=1}^{n}p_{k}|u_{k}|\\
&\leq\bigvee\limits_{l_{k,j}=1}^{2^{m}}\sum\limits_{k=1}^{n}\prod\limits_{j=1}^{p_{k}}t_{k,j}(\frac{l_{k,j}\pi}{2^{m}})u_{k}+\frac{\pi}{2^{m}}\sum\limits_{k=1}^{n}p_{k}|u_{k}|.
\end{align*}
\noindent
Since this is true for every $\theta_{k,j}\in[0,\frac{\pi}{2}]\ (k\in\lbrace 1,...,n\rbrace, j\in\lbrace 1,...,p_{k}\rbrace)$ we have
\begin{equation*}
0\leq f-\bigvee\limits_{l_{k,j}=1}^{2^{m}}\sum\limits_{k=1}^{n}\prod\limits_{j=1}^{p_{k}}t_{k}(\frac{l_{k,j}\pi}{2^{m}})\leq\frac{\pi}{2^{m}}\sum\limits_{k=1}^{n}p_{k}|u_{k}|.
\end{equation*}
\noindent
It follows that the sequence $\sigma_{m}:=\bigvee\limits_{l_{k,j}=1}^{2^{m}}\sum\limits_{k=1}^{n}\prod\limits_{j=1}^{p_{k}}t_{k}(\frac{l_{k,j}\pi}{2^{m}})$ converges relatively uniformly to $f$. Finally, for $f\in E$, there exist sequences $(a_{n}),(b_{n})$ in $E$ such that $a_{n}\overset{ru}{\rightarrow}f^{+}$ and $b_{n}\overset{ru}{\rightarrow}f^{-}$. Then $a_{n}-b_{n}\overset{ru}{\rightarrow}f$.
\end{proof}

\smallskip
We are ready to deal with the Archimedean tensor product of Archimedean vector lattices over $\mathbb{K}$. Parts (1), (2), and (4) of the following theorem extend the corresponding parts of Theorem 4.2 in \cite{Fr} and Lemma 4.2. Parts (2) and (4) generalize corresponding results by Schep for real Archimedean vector lattices in Section 2 of \cite{Sh}. Part (3) is slightly weaker than the real analogues found in \cite{Fr} and \cite{Sh}, but it sufficient for obtaining our results for maps of order bounded variation in the next section. We do not yet know if it is possible to strengthen (3).

\smallskip
\noindent
\th{Theorem 4.7.} Let $E_{1},...,E_{s}$ be
Archimedean vector lattices over $\mathbb{K}$.
\begin{enumerate}
\item[(1)] \textit{There exists an essentially unique Archimedean
vector lattice $\bar{\otimes }_{k=1}^{s}E_{k}$\textit{\ over $\mathbb{K}$ and a vector
lattice }$s$\textit{-morphism }$\bar{\otimes }:\times
_{k=1}^{s}E_{k}\rightarrow \bar{\otimes }_{k=1}^{s}E_{k}$\textit{\ such
that for every Archimedean vector lattice }$F$\textit{\ over }$\mathbb{K}$ 
\textit{and every vector lattice }$s$\textit{-morphism }$T:\times
_{k=1}^{s}E_{k}\rightarrow F$\textit{, there exists a unique vector lattice
homomorphism }$T^{\bar{\otimes }}:\bar{\otimes }
_{k=1}^{s}E_{k}\rightarrow F$\textit{\ such that }$T^{\bar{\otimes }
}\circ \bar{\otimes }=T$\textit{. } }

\item[(2)] \textit{There exists an injective $\mathbb{K}$-linear map $S:\otimes
_{k=1}^{s}E_{k}\rightarrow $\textit{\ }$\bar{\otimes }_{k=1}^{s}E_{k}$
such that $S\circ\otimes=\bar{\otimes}$. }

\item[(3)] $\tau(\otimes_{k=1}^{s}E_{k},\bar{\otimes}_{k=1}^{s}E_{k})\leq 2$. Thus, $\otimes_{k=1}^{s}E_{k}$ is dense in $\bar{\otimes}_{k=1}^{s}E_{k}$ in the relatively uniform topology.

\item[(4)] \textit{For every} $w\in (\bar{\otimes }_{k=1}^{s}E_{k})\setminus \{0\}$\textit{\ there exist } $x_{1}\otimes\dots\otimes x_{s}\in\otimes_{k=1}^{s} E_{k\rho}$ with $x_{k}\in E_{k}^{+}\ (k\in
\{1,...,s\})$\textit{\ such that }$0<(x_{1}\otimes \dots \otimes x_{s})\leq
|w|$\textit{, i.e. }$\otimes _{k=1}^{s}E_{k}$\textit{\ is order dense in }$\bar{\otimes }_{k=1}^{s}E_{k}$\textit{.}
\end{enumerate}
\endth

\begin{proof} By Lemma 4.2, statements (1)-(4) are valid for $\mathbb{K}=\mathbb{R}$. We assume in the proof below that $\mathbb{K}=\mathbb{C}$.

(1) Let $E_{1},...,E_{s},F$ be Archimedean vector lattices over $\mathbb{C}$. Denote by $(\bar{\otimes}_{k=1}^{s}E_{k\rho },\bar{\otimes })$
the Archimedean vector lattice tensor product of $E_{1\rho},...,E_{s\rho}$. We
claim that the pair $((\bar{\otimes }_{k=1}^{s}E_{k\rho })_{|\mathbb{C}|},\bar{\otimes }_{|\mathbb{C}|})$ is the unique Archimedean complex vector lattice
tensor product of $E_{1},...,E_{s}$. Let $T:\times
_{k=1}^{s}E_{k}\rightarrow F$ be a vector lattice $s$-morphism. From Lemma
4.2, the map $\bar{\otimes}$ induces a unique vector lattice
homomorphism $T_{\rho}^{\bar{\otimes}}$ on $\bar{\otimes}_{k=1}^{s}E_{k\rho}$ such that $T_{\rho}^{\bar{\otimes}}\circ 
\bar{\otimes}=T_{\rho}$. Also, the map $T_{\rho}^{\bar{\otimes}}$ extends uniquely to a vector lattice homomorphism $(T_{\rho}^{\bar{\otimes}})^{\mu_{2,4}}$ on $(\bar{\otimes}_{k=1}^{s}E_{k\rho})^{\mu_{2,4}}$ (Corollary 3.1). By Proposition 3.4(1), the map $\bar{\otimes}_{|\mathbb{C}|}$ is a vector lattice $s$-morphism and $(T_{\rho}^{\bar{\otimes }})_{|\mathbb{C}|}$ is a vector lattice
homomorphism. We will prove that the map $(T_{\rho}^{\bar{\otimes}})_{|\mathbb{C}|}:(\bar{\otimes}_{k=1}^{s}E_{k\rho})_{|\mathbb{C}|}\rightarrow F$ is the unique vector lattice homomorphism such that $(T_{\rho}^{\bar{\otimes}})_{|\mathbb{C}|}\circ\bar{\otimes}_{|\mathbb{C}|}=T$. Indeed, for every $(f_{0}^{1}+if_{1}^{1},...,f_{0}^{s}+if_{1}^{s})\in\times
_{k=1}^{s}E_{k}$ we have
\begin{eqnarray*}
(T_{\rho}^{\bar{\otimes }})_{|\mathbb{C}|}\circ \bar{\otimes }_{|%
\mathbb{C}|}(f_{0}^{1}+if_{1}^{1},...,f_{0}^{s}+if_{1}^{s}) &=&(T_{\rho }^{%
\bar{\otimes }})_{|\mathbb{C}|}(\sum\limits_{\epsilon _{k}\in \{0,1\}}%
\bar{\otimes }(f_{\epsilon _{1}}^{1},...,f_{\epsilon
_{s}}^{s})i^{\sum\limits_{k=1}^{s}\epsilon _{k}}) \\
&=&\sum\limits_{\epsilon _{k}\in \{0,1\}}T_{\rho }^{\bar{\otimes }%
}\circ \bar{\otimes }(f_{\epsilon _{1}}^{1},...,f_{\epsilon
_{s}}^{s})i^{\sum\limits_{k=1}^{s}\epsilon _{k}} \\
&=&\sum\limits_{\epsilon _{k}\in \{0,1\}}T_{\rho }(f_{\epsilon
_{1}}^{1},...,f_{\epsilon _{s}}^{s})i^{\sum\limits_{k=1}^{s}\epsilon _{k}} \\
&=&T(f_{0}^{1}+if_{1}^{1},...,f_{0}^{s}+if_{1}^{s}).
\end{eqnarray*}
\noindent
Since every vector lattice $\mathbb{C}$-homomorphism is
real, the uniqueness of $(T_{\rho }^{\bar{\otimes}})_{|\mathbb{C}|}$
follows from the uniqueness of $(T_{\rho }^{\bar{\otimes}})^{\mu_{2,4}}$.

The proof of uniqueness of the Archimedean complex vector lattice tensor product is not different from the real case.

(2) Consider the newly minted tensor product $(\bar{\otimes }_{k=1}^{s}E_{k},\bar{\otimes })$ constructed in (1). By Lemma 4.2, there exists an
Archimedean vector lattice $G$ over $\mathbb{R}$ and a vector lattice $s$-morphism $T:\times_{k=1}^{s}E_{k\rho}\rightarrow G$ such that the induced
linear map $T^{\otimes}:\otimes_{k=1}^{s}E_{k\rho}\rightarrow G$ is injective. By taking the square mean completion of $G$, if necessary, we will assume that $G$ is square mean complete. By taking vector space complexifications, we find an injective vector lattice $s$-morphism $(T^{\otimes})_{\mathbb{C}}:(\otimes_{k=1}^{s}E_{k\rho})_{\mathbb{C}}\rightarrow G_{\mathbb{C}}$, or equivalently by Lemma 4.3, $(T^{\otimes})_{\mathbb{C}}:\otimes_{k=1}^{s}E_{k}\rightarrow G_{\mathbb{C}}$. Moreover, if $(T_{\mathbb{C}})^{\otimes}:\otimes_{k=1}^{s}E_{k}\rightarrow G$ is the unique linear map induced by $T_{\mathbb{C}}$, then for $f_{0}^{k}+if_{1}^{k}\in E_{k}\ (k\in\lbrace 1,...,s\rbrace)$,
\begin{align*}
(T_{\mathbb{C}})^{\otimes}(f_{0}^{1}+if_{1}^{1}\otimes\dots\otimes f_{0}^{s}+if_{1}^{s})&=T_{\mathbb{C}}(f_{0}^{1}+if_{1}^{1},...,f_{0}^{s}+if_{1}^{s})\\
&=\sum\limits_{\epsilon _{k}\in \{0,1\}}T(f_{\epsilon_{1}}^{1},...,f_{\epsilon _{s}}^{s})i^{\sum\limits_{k=1}^{s}\epsilon_{k}}\\
&=\sum\limits_{\epsilon _{k}\in \{0,1\}}T^{\otimes}(f_{\epsilon_{1}}^{1}\otimes\dots\otimes f_{\epsilon _{s}}^{s})i^{\sum\limits_{k=1}^{s}\epsilon_{k}}.
\end{align*}
\noindent
In particular, $((T_{\mathbb{C}})^{\otimes})$ is a real map with $((T_{\mathbb{C}})^{\otimes})_{\rho}=T^{\otimes}$, and therefore we have $(T_{\mathbb{C}})^{\otimes}=(T^{\otimes})_{\mathbb{C}}$. From part (1) of this theorem there exists a unique vector lattice $\mathbb{C}$-homomorphism $(T_{\mathbb{C}})^{\bar{\otimes}}:\bar{\otimes }_{k=1}^{s}E_{k}\rightarrow G_{\mathbb{C}}$ such that $(T_{\mathbb{C}})^{\bar{\otimes }}\circ \bar{\otimes}=T_{\mathbb{C}}$.
Moreover, there exists a unique $\mathbb{C}$-linear map $S:\otimes_{k=1}^{s}E_{k}\rightarrow\bar{\otimes }_{k=1}^{s}E_{k}$
such that $S\circ\otimes=\bar{\otimes}$. Then $(T_{\mathbb{C}})^{\bar{\otimes}}\circ S\circ\otimes=T_{\mathbb{C}}$, and hence $(T_{\mathbb{C}})^{\bar{\otimes}}\circ S=(T_{\mathbb{C}})^{\otimes}=(T^{\otimes})_{\mathbb{C}}$.
Therefore $S$ is injective.

(3) By Lemma 4.2 we know that $\otimes_{k=1}^{s}E_{k\rho}$ is relatively uniformly dense in $\bar{\otimes}_{k=1}^{s}E_{k\rho}$. We also know from Proposition 4.6 that $\bar{\otimes}_{k=1}^{s}E_{k\rho}$ is relatively uniformly dense in $(\bar{\otimes}_{k=1}^{s}E_{k\rho})^{\mu_{2,4}}$. By taking vector space complexifications, we have $\tau(\otimes_{k=1}^{s}E_{k},\bar{\otimes}_{k=1}^{s}E_{k})\leq 2$.

(4) Suppose $w\in(\bar{\otimes}_{k=1}^{s}E_{k})\setminus\{0\}$. Then $0<|w|\in(\bar{\otimes}_{k=1}^{s}E_{k\rho})^{\mu_{2,4}}$. Since $\bar{\otimes}_{k=1}^{s}E_{k\rho}$ is order dense in $(\bar{\otimes}_{k=1}^{s}E_{k\rho})^{\delta}$, it is also order dense in $(\bar{\otimes}_{k=1}^{s}E_{k\rho})^{\mu_{2,4}}$. Thus there exists $w_{0}\in\bar{\otimes}_{k=1}^{s}E_{k\rho}$ such that $0<w_{0}\leq|w|$. From Lemma 4.2, there exists $x_{1}\otimes\dots\otimes x_{s}\in\otimes_{k=1}^{s}E_{k\rho}$ with $x_{k}\in E_{k}^{+}\ (k\in \{1,...,s\})$ such that $0<(x_{1}\otimes\dots\otimes x_{s})\leq w_{0}$.
\end{proof}

\smallskip In (1) it is necessary to take the vector lattice
complexification of $\bar{\otimes}_{k=1}^{s}E_{k\rho }$ to ensure that 
$(\bar{\otimes}_{k=1}^{s}E_{k\rho})_{|\mathbb{C}|}$ is an Archimedean vector
lattice over $\mathbb{C}$. Indeed, Theorems 4.10 and 4.11 furnish examples where the
vector space complexification $\left(\bar{\otimes}_{k=1}^{s}E_{k\rho}\right)_{\mathbb{C}}$ does not suffice. We need two lemmas first.

\smallskip
\noindent
\th{Lemma 4.8.} Let $X$ and $Y$ be nonempty
subsets of $\mathbb{R}$ without isolated points. Then the function $%
S:(x,y)\mapsto \sqrt{x^{2}+y^{2}}\ ((x,y)\in X\times Y)$ is in the square
mean completion of $C(X)\bar{\otimes }C(Y)$ but for all nonempty open subsets $U$ of $X$ and $W$ of $Y$ we have $S\mid _{U\times
W}\notin C(U)\otimes C(W)$.
\endth

\begin{proof} For $f\in C(X)$ and $g\in C(Y)$ we identify $f\otimes g$ with the
function $(x,y)\mapsto f(x)g(y)$ $((x,y)\in X\times Y)$. Consider the
element $S$ of the square mean completion of $C(X)\bar{\otimes}C(Y)$ defined
by
\begin{equation*}
(x,y)\mapsto \sqrt{x^{2}+y^{2}}\ ((x,y)\in X\times Y)\text{.}
\end{equation*}
\noindent
Let $U$ and $W$ be open nonempty subsets of $X$ and $Y$, respectively. We will show that the vector subspace of $C(U)$ generated
by $\{S(\cdot ,y):y\in W\}$, whose elements are considered as functions on $U$, is not finite-dimensional. It follows (see Proposition 1 in \cite{Hag})
that $S\mid _{U\times W}\notin C(U)\otimes C(W)$. Since $W$ is open and
nonempty and $Y$ has no isolated points, we can choose $\alpha _{k}\in W$
(for all $k\in \mathbb{N}$) for which $\alpha _{i}^{2}\neq \alpha _{j}^{2}$
when $i\neq j$. Let $n\in \mathbb{N}$ and let $\lambda _{1},...,\lambda
_{n}\in \mathbb{R}$ for which $\lambda _{k}\sqrt{x^{2}+\alpha _{k}^{2}}=$ $\lambda _{k}S(x,y_{k})=0$ for all $x\in U$. Since the function $x\mapsto
\lambda _{k}\sqrt{x^{2}+\alpha _{k}^{2}}\text{ }(x\in \mathbb{R})$ is $n$
times differentiable at every $x\in X$, a routine calculation shows that the $n\times n$ matrix $A(x)$ defined by $A(x)_{ij}=
\frac{1}{(x^{2}+\alpha _{j}^{2})^{\frac{2i-1}{2}}}$ when evaluated at the
vector $(\lambda _{1},...,\lambda _{n})$ yields the vector $(0,...,0)$ for
every non-zero $x\in U$. However, $\prod\limits_{k=1}^{n}\sqrt{x^{2}+\alpha _{k}^{2}}
\det (A(x))=\det (B(x))$ where the $n\times n$ matrix $B(x)$ is defined by $
B(x)_{ij}=\frac{1}{(x^{2}+\alpha _{j}^{2})^{i-1}}$, which has (Vandermonde)
determinant
\begin{equation*}
\prod\limits_{1\leq j<k\leq n}\Bigl(\frac{1}{x^{2}+\alpha_{j}^{2}}-\frac{1}{x^{2}+\alpha_{k}^{2}}\Bigr)=\prod\limits_{1\leq j<k\leq n}\frac{\alpha
_{j}^{2}-\alpha _{k}^{2}}{(x^{2}+\alpha _{j}^{2})(x^{2}+\alpha _{k}^{2})}\neq 0\text{.}
\end{equation*}
\noindent
Thus $\det (A(x))\neq 0$ for every non-zero $x\in U$, the vector subspace of $C(U)$ generated by $\{S(\cdot ,y):y\in Y\}$ (as functions on $U$) is
infinite dimensional, and $S\mid _{U\times W}\notin C(U)\otimes C(W)$.
\end{proof}

\smallskip
\noindent
\th{Lemma 4.9.} Let $X$ and $Y$ be
nonempty subsets of $\mathbb{R}$ without isolated points and let $f\in C(X)\bar{\otimes }C(Y)$. Then there exists a nonempty open subset $V$ of $X\times Y$ and $g\in C(X)\otimes C(Y)$ such that $f|_{V}=g|_{V}$.
\endth

\begin{proof}
Note that $C(X)\bar{\otimes }C(Y)$ is the vector lattice
generated by $C(X)\otimes C(Y)$ in $C(X\times Y)$\ (\cite{Fr}, Section 4).
Every element $f\in C(X)\bar{\otimes }C(Y)$ is of the form $%
f=\bigwedge\limits_{j=1}^{n}\bigvee\limits_{k=1}^{m}f_{j,k}$ where $%
f_{j,k}\in C(X)\otimes C(Y)$ for each $j$ and each $k$ (see Exercise 4.8 in 
\cite{AB}). Let $f_{1},f_{2}\in C(X)\otimes C(Y)$. If $f_{1}\neq f_{2}
$, we may assume there exists $(x,y)$ such that $f_{1}(x,y)<f_{2}(x,y)$ and
then there exists a nonempty open subset $O$ of $X\times Y$ such that $f_{1}\wedge f_{2}=f_{1}$ on $O$. Of course such an open subset $O$ also
exists if $f_{1}=f_{2}$. By repeating this argument there exists a nonempty
open set $U\subseteq O$ such that $\bigwedge\limits_{j=1}^{n}(\bigvee\limits_{k=1}^{m}f_{j,k})=\bigvee\limits_{k=1}^{m}f_{j_{0},k}$ on $U$ for some $j_{0}\in\lbrace 1,...,n\rbrace$. Similarly, there exists a nonempty open set $V\subset U$ such that $\bigvee\limits_{k=1}^{m}f_{j_{0},k_{0}}=f_{j_{0},k_{0}}$ on $V$ for some $k_{0}\in\lbrace 1,...,m\rbrace$.
\end{proof}

\smallskip
\noindent
\th{Theorem 4.10.} Let $X$ and $Y$ be
nonempty subsets of $\mathbb{R}$ without isolated points. Then $C(X)\bar{\otimes}C(Y)$ is not square mean complete. Therefore $(C(X)\bar{\otimes}C(Y))_{\mathbb{C}}$ is not an Archimedean vector lattice over $\mathbb{C}$.
\endth

\begin{proof} Assume that the element $S$ of Lemma 4.8 is in $C(X)\bar{\otimes}C(Y)$. Then by Lemma 4.9 there exists a nonempty open set $V$ in $X\times Y$ and an element $g\in C(X)\otimes C(Y)$ such that $g\mid_{V}=S\mid_{V}.$ However the open set $V$ contains a nonempty open subset
of the form $U\times W$ with $0\notin U$. This contradicts Lemma 4.8.
\end{proof}

\smallskip
We use Theorem 4.10 to prove the following.

\smallskip
\noindent
\th{Theorem 4.11.} If $X$ and $Y$ are uncountable compact metrizable spaces then $C(X)\bar{\otimes}C(Y)$ is not square mean complete. Therefore $(C(X)\bar{\otimes}C(Y))_{\mathbb{C}}$ is not an Archimedean vector lattice over $\mathbb{C}$.
\endth

\begin{proof} By Theorem 1 in \cite{Pel}, we know that
both $X$ and $Y$ contain a closed subset homeomorphic with the Cantor set $\mathbb{D}$. Then $\mathbb{D\times D}$ can be viewed as a closed subset of $X\times Y$ and the function $F_{0}:(x,y)\longmapsto \sqrt{x^{2}+y^{2}}$($(x,y)\in \mathbb{D\times D}$) is continuous. By Tietze's Extension Theorem, the function $x\longmapsto x$\ ($x\in \mathbb{D}$) can be extended to continuous functions $f$ and $g$ on $X$ and $Y$, respectively. Then the function $F:(x,y)\longmapsto \sqrt{f(x)^{2}+g(y)^{2}}\ ((x,y)\in X\mathbb{\times}Y)$
is a continuous function in the square mean completion of $C(X)\bar{\otimes}C(Y)$ that extends $F_{0}$. If $F$ were in $C(X)\bar{\otimes}C(Y)$ itself
then its restriction to $\mathbb{D\times D}$ would be in $C(\mathbb{D})\bar{\otimes}C(\mathbb{D})$ which by Theorem 4.10 is impossible. This proves the
theorem.
\end{proof}

\smallskip
It is certainly tempting to conjecture the following.

\smallskip
\noindent
\th{Conjecture 4.12.} If $X$ and $Y$ are infinite compact metrizable spaces then $C(X)\bar{\otimes}C(Y)$ is not square mean complete.
\endth

\smallskip
The above two theorems show that the old way of complexifying
Archimedean vector lattices via vector space complexifications is inadequate for pursuing complex analysis on Archimedean complex vector lattices.

We remark that the complex Archimedean vector lattice
tensor product, like its real counterpart (\cite{Fr}, Theorem 5.3, \cite{Sh}, Section 2), possesses as well a universal property with respect to
positive multilinear maps and complex uniformly complete vector lattices as
range. The proof of this universal property, stated in the theorem below, is similar to the proof of Theorem 4.7(1) and is left to reader.

\smallskip
\noindent
\th{Theorem 4.13.} Let $E_{1},...,E_{s},F$ be Archimedean vector lattices over $\mathbb{K}$ with $F$ uniformly complete. If $T:\times _{k=1}^{s}E_{k}\rightarrow F$ is a positive $s_{\mathbb{K}}$-linear map, then there exists a unique positive $\mathbb{K}$-linear map $T^{\bar{\otimes }}:\bar{\otimes }_{k=1}^{s}E_{k}\rightarrow F$ such that $T^{\bar{\otimes }}\circ\bar{\otimes}=T$.
\endth

\smallskip
A reformulation of part (1) of Theorem 4.7 in terms of Archimedean real vector lattices and vector lattice complexifications is the following.

\smallskip
\noindent
\th{Theorem 4.14.} Let $E_{1},...E_{s},F$ be Archimedean vector lattices over $\mathbb{R}$ and suppose that $T:\times_{k=1}^{s}E_{k}\rightarrow F$ is a
vector lattice $s_{\mathbb{R}}$-morphism. There exists a unique vector
lattice $s_{\mathbb{C}}$-morphism $(T_{|\mathbb{C}|})^{\bar{\otimes}}:\bar{\otimes}_{k=1}^{s}E_{k|\mathbb{C}|}\rightarrow F_{|\mathbb{C}|}$
such that $(T_{|\mathbb{C}|})^{\bar{\otimes}}\circ\bar{\otimes}|_{\times_{k=1}^{s}E_{k}}=T$.
\endth

\begin{proof} Consider $T$ to be a vector lattice $s_{\mathbb{R}}$-morphism from $\times_{k=1}^{s}E_{k}$ to $F^{\mu_{2,4}}$.  By Proposition 3.4(1) there 
exists a unique vector lattice $s_{\mathbb{C}}$-morphism $T_{\left\vert\mathbb{C}\right\vert}:\times_{k=1}^{s}E_{k\left\vert \mathbb{C}\right\vert}\rightarrow F_{\left\vert\mathbb{C}\right\vert}$ such that $T_{|\mathbb{C}|}|_{\times_{k=1}^{s}E_{k}}=T$. If $(T_{|\mathbb{C}|})^{\bar{\otimes}}$ is the unique vector lattice $\mathbb{C}$-homomorphism induced by $T_{\mathbb{C}}$ then $(T_{|\mathbb{C}|})^{\bar{\otimes}}\circ\bar{\otimes}=T_{|\mathbb{C}|}$. In particular, $(T_{|\mathbb{C}|})^{\bar{\otimes}}\circ\bar{\otimes}|_{\times_{k=1}^{s}E_{k}}=T$.
\end{proof}

\vspace{.5in}
\noindent
\textbf{5. Applications of the Archimedean Vector Lattice Tensor Product}
\vspace{.2in}

In this section, we give some applications of the Archimedean vector lattice tensor product of Archimedean vector lattices over $\mathbb{K}$. Indeed, we use this tensor product to prove the existence of $s$-powers for every $s\in \mathbb{N}\backslash \{1\}$ and every Archimedean vector lattice over $\mathbb{K}$. We also generalize Theorem 3.1 of \cite{BusvR} to $s_{\mathbb{K}}$-linear maps of order bounded variation.

\smallskip
Central to the theory of $s$-powers are orthosymmetric $s$-morphisms.

\smallskip
\noindent
\th{Definition 5.1.} For Archimedean vector lattices 
$E_{1},...,E_{s},F$ over $\mathbb{K}$, a map $T:\times _{s}E\rightarrow F$ is called orthosymmetric if $T(f_{1},...,f_{s})=0$ whenever there exist $i,j\in \{1,...,s\}$ such that $%
|f_{i}|\wedge |f_{j}|=0$.
\endth

\smallskip
\noindent
\th{Definition 5.2.} Let $E$ be an Archimedean
vector lattice over $\mathbb{K}$ and let $s\in \mathbb{N}\setminus \{1\}$.
We call a pair $(E^{\circledS },\circledS )$ an $s$-power of $E$ if the following hold.
\begin{enumerate}
\item[(1)] \textit{$E^{\circledS }$\textit{\ is an Archimedean vector
lattice over }$\mathbb{K}$\textit{. } }

\item[(2)] \textit{$\circledS:\times_{s}E\rightarrow E^{\circledS }$\textit{\ is an orthosymmetric vector lattice }$s$\textit{-morphism.} }

\item[(3)] \textit{For every Archimedean vector lattice $F$\textit{\ over }$
\mathbb{K}$\textit{\ and every orthosymmetric }$s$\textit{-morphism }$
T:\times _{s}E\rightarrow F$\textit{, there exists a unique vector lattice
homomorphism }$T^{\circledS }$\textit{\ such that }$T^{\circledS }\circ
\circledS =T$\textit{.} }
\end{enumerate}
\endth

\smallskip
We address the existence and uniqueness of $s$-powers for
Archimedean vector lattices over $\mathbb{K}$ in our next theorem, which
extends Theorem 3.2 in \cite{BoBus}. We denote $E\bar{\otimes}\dots\bar{\otimes}E\ (s\ \text{times})$ by $\bar{\otimes}_{s}E$.

\smallskip
\noindent
\th{Theorem 5.3.} If $E$ is an
Archimedean vector lattice over $\mathbb{K}$ and $s\in\mathbb{N}\setminus\lbrace 1\rbrace$ then there exists an $s$-power
of $E$, unique up to vector lattice isomorphism.
\endth

\begin{proof}
The proof for $\mathbb{K}=\mathbb{R}$ is the content of the proof of Theorem 3.2 in \cite{BoBus}, so we assume $\mathbb{K}=\mathbb{C}$. Let $E$ be an Archimedean vector lattice over $\mathbb{C}$ and let $I$ be the smallest uniformly closed ideal of $\bar{\otimes}_{s}E$ that contains
\begin{equation*}
\bigl\lbrace f_{1}\bar{\otimes}\dots\bar{\otimes}f_{s}:f_{1},...,f_{s}\in E\ \text{and}\ |f_{j}|\wedge|f_{k}|=0\ \text{for some}\ j,k\in\lbrace 1,...,s\rbrace\bigr\rbrace.
\end{equation*}
\noindent
Given $f\in\bar{\otimes}_{s}E$, we denote the equivalence class of $f$ in $(\bar{\otimes}_{s}E)/I$ by $[f]$. Then $(\bar{\otimes}_{s}E)/I$ is a vector space over $\mathbb{C}$ under the operations $[f]+[g]=[f+g]$ and $\lambda[f]=[\lambda f]$ for all $\lambda\in\mathbb{C}$. From $[f+ig]=[f]+i[g]$, we see that $(\bar{\otimes}E_{s})/I=\bigl((\bar{\otimes}E_{s})_{\rho}/I_{\rho}\bigr)_{\mathbb{C}}$ (see page 198 of \cite{Zan2}). Also, $I_{\rho}$ is a uniformly closed ideal in $(\bar{\otimes}E_{s})_{\rho}$, and thus $(\bar{\otimes}E_{s})_{\rho}/I_{\rho}$ is an Archimedean vector lattice over $\mathbb{R}$. Let $p:(\bar{\otimes}E_{s})_{\rho}\rightarrow(\bar{\otimes}E_{s})_{\rho}/I_{\rho}$ be the natural vector lattice homomorphism, i.e., $p(f)=[f]$ for all $f\in(\bar{\otimes}E_{s})_{\rho}$. Let $[f],[g]\in(\bar{\otimes}E_{s})_{\rho}/I_{\rho}$. By Theorem 3.12 in \cite{BusSch}, $\mu_{2,4}([f],[g])=[\mu_{2,4}(f,g)]$, which is in $(\bar{\otimes}E_{s})_{\rho}/I_{\rho}$ since $(\bar{\otimes}E_{s})_{\rho}$ is square mean complete. Hence $(\bar{\otimes}E_{s})_{\rho}/I_{\rho}$ is also square mean complete and $(\bar{\otimes}E_{s})/I$ is an Archimedean complex vector lattice. Next, let $q:\bar{\otimes}E_{s}\rightarrow(\bar{\otimes}E_{s})/I$ be the natural vector lattice homomorphism from $\bar{\otimes}E_{s}$ to $(\bar{\otimes}E_{s})/I$. Following Theorem 4 in \cite{BusvR2}, it is straightforward to show that $\bigl((\bar{\otimes}E_{s})/I,q\circ\bar{\otimes}\bigr)$ is an $s$-power of $E$. The proof of the uniqueness of $s$-powers is the same as the real case.
\end{proof}

\smallskip
The Archimedean complex vector lattice tensor product can also be used to
obtain results for multilinear maps over $\mathbb{K}$ of ordered bounded
variation. We start with some definitions.

Let $E$ be an Archimedean vector lattice over $\mathbb{K}$ and let $a\in E^{+}$. A \textit{partition} of $a$ is a finite sequence $%
\{x_{k}\}_{k=1}^{n}$ in $E^{+}$ such that $\sum\limits_{k=1}^{n}x_{k}=a$. As
in \cite{BusvR}, we denote the set of all partitions of $a$ by $\prod a$ and
abbreviate a partition $\{x_{k}\}_{k=1}^n$ of $a$ by $x$, which explains the
shorthand $x\in \prod a$.

\smallskip
\noindent
\th{Definition 5.4.} Let $E_{1},...,E_{s},F$ be Archimedean vector lattices over $\mathbb{K}$. We say
that an $s_{\mathbb{K}}$-linear map $T:\times _{k=1}^{s}E_{k}\rightarrow F$
is of order bounded variation if for all $a_{k}\in E_{k}^{+}\ (k\in
\{1,...,s\})$ the set
\begin{align*}
\bigl\lbrace\sum\limits_{n_{1},...,n_{s}}|T(x^{1}_{n_{1}},..,x^{s}_{n_{s}})|:x^{k}\in\prod a_{k} (k=1,...,s)\rbrace
\end{align*}

\noindent is order bounded. We denote by $\mathcal{L}_{bv}(E_{1},...,E_{s};F)$ the space of all $s_{\mathbb{K}}$-linear maps of
order bounded variation from $\times _{k=1}^{s}E_{k}$ into $F$.
\endth

\smallskip
\noindent
\th{Definition 5.5.} Let $V$ be a vector
space over $\mathbb{K}$. We call $K\subseteq V$ a cone in $V$ if $K+K\subseteq K$, $\lambda K\subseteq K$ for every $\lambda\in\mathbb{K}^{+}$, and $K\cap(-K)=\lbrace 0\rbrace$. The pair $(V,K)$ where $V$ is a vector space over $\mathbb{K}$ and $K$ is a cone in $V$ is called an ordered vector space over $\mathbb{K}$.
\endth

\smallskip
For example, if $E_{1},...,E_{s},F$ are Archimedean
vector lattices over $\mathbb{K}$ then $\mathcal{L}_{bv}(E_{1},...,E_{s};F)$
is a vector space over $\mathbb{K}$ and the set of all positive maps in $\mathcal{L}_{bv}(E_{1},...,E_{s};F)$, which we denote by $\mathcal{L}^{+}_{bv}(E_{1},...,E_{s};F)$, is a cone.

Let $\left(V_{1},K_{1}\right) ,...,\left(V_{s},K_{s}\right),(W,K)$ be ordered vector spaces over $\mathbb{K}$. We say that a map $T:\times
_{k=1}^{s}V_{k}\rightarrow W$ is \textit{positive} if $T(\times _{k=1}^{s}K_{k})\subseteq K$. If a positive $\mathbb{K}$-linear map $T$ is bijective and has a positive inverse, we call $T$ an 
\textit{ordered vector space isomorphism}. If there exists an ordered vector
space isomorphism between ordered vector spaces $(V,K)$ and $(W,K^{\prime})$
over $\mathbb{K}$ we say that $(V,K)$ and $(W,K^{\prime})$ are \textit{isomorphic as ordered vector spaces}.

For the proof of the following lemma, let $E$ be an Archimedean vector lattice over $\mathbb{K}$, let $(V,K)$ be an ordered vector space over $\mathbb{K}$, and let $\phi:E\rightarrow V$ be an ordered vector space isomorphism with respect to the cones $E^{+}$ and $K$. It is readily checked that the map $m:V\rightarrow V$ defined by $m(v)=\phi(|\phi^{-1}(v)|)$ is an Archimedean modulus on $V$ with $m(V)=K$.

\smallskip
\noindent
\th{Lemma 5.6.} If $E$ is an Archimedean vector lattices over $\mathbb{K}$, $(V,K)$ is an ordered vector space over $\mathbb{K}$, and $\phi:E\rightarrow V$ is an ordered vector space isomorphism with respect to the cones $E^{+}$ and $K$ then

\begin{itemize}
\item[(1)] $V$ is an Archimedean vector lattice over $\mathbb{K}$ with $K$ as positive cone, and
\item[(2)] $\phi$ is a vector lattice isomorphism.
\end{itemize}
\endth

\smallskip
The following result generalizes Proposition 3.2(4) in \cite{dS} as well as Theorem 3.1 in \cite{BusvR}.

\smallskip
\noindent
\th{Theorem 5.7.} Let $E_{1},...,E_{s},F$ be Archimedean vector lattices over $\mathbb{K}$ with $F$ Dedekind
complete.
\begin{enumerate}
\item[(1)] For any $s_{\mathbb{K}}$-linear map of order bounded
variation $T:\times_{k=1}^{s}E_{k}\rightarrow F$ there exists a unique
order bounded $\mathbb{K}$-linear map $T^{\bar{\otimes }}:\bar{\otimes }_{k=1}^{s}E_{k}\rightarrow F$ such that $T(f_{1},..,f_{s})=T^{\bar{\otimes }}(f_{1}\bar{\otimes }\dots \bar{\otimes }f_{s})$
for every $f_{k}\in E_{k}\ (k\in \{1,...,s\})$.

\item[(2)] $\mathcal{L}_{bv}(E_{1},...,E_{s};F)$ is a Dedekind complete Archimedean
vector lattice over $\mathbb{K}$ and the correspondence $T\mapsto T^{\bar{\otimes }}$ is a vector lattice
isomorphism from $\mathcal{L}_{bv}(E_{1},...,E_{s};F)$ onto $\mathcal{L}_{b}(\bar{\otimes }_{k=1}^{s}E_{k},F)$.

\item[(3)] For $T\in\mathcal{L}_{bv}(E_{1},...,E_{s};F)$,

$|T|(a_{1},...,a_{n})=\sup\lbrace\sum\limits_{n_{1},...,n_{s}}|T(x_{n_{1}}^{1},...,x_{n_{s}}^{s})|:x^{k}\in
\prod a_{k}\ (k\in\lbrace 1,...,s)\rbrace$\newline
\noindent
for all $a_{k}\in E_{k}^{+}\ (k\in \{1,...,s\})$.
\end{enumerate}
\endth

\begin{proof}
For $\mathbb{K}=\mathbb{R}$ the result is Theorem 3.1 in \cite{BusvR}. We thus assume $\mathbb{K}=\mathbb{C}$.

(1) For the uniqueness, suppose that $T:\times_{k=1}^{s}E_{k}\rightarrow F$ is an $s_{\mathbb{C}}$-linear map of order bounded variation, and assume that $S_{1}, S_{2}$ are complex order bounded linear maps from $\bar{\otimes}_{k=1}^{s}E_{k}$ to $F$ such that $T=S_{1}\circ\bar{\otimes}$ and $T=S_{2}\circ\bar{\otimes}$. Then $S_{1}-S_{2}=0$ identically on $\otimes_{k=1}^{s}E_{k}$. By relatively uniform density, $S_{1}-S_{2}=0$ on $(\otimes_{k=1}^{s}E_{k})_{2}$, where $(\otimes_{k=1}^{s}E_{k})_{2}$ denotes the pseudo uniform closure of $\otimes_{k=1}^{s}E_{k}$ in $\bar{\otimes}_{k=1}^{s}E_{k}$ (see Section 2). By relatively uniform density again, $S_{1}-S_{2}=0$ on $(\otimes_{k=1}^{s}E_{k})_{3}$, where $(\otimes_{k=1}^{s}E_{k})_{3}$ denotes the pseudo uniform closure of $(\otimes_{k=1}^{s}E_{k})_{2}$ in $\bar{\otimes}_{k=1}^{s}E_{k}$. From Theorem 4.7(3), we have $(\otimes_{k=1}^{s}E_{k})_{3}=\bar{\otimes}_{k=1}^{s}E_{k}$. We next turn to the existence. Define $\bar{T}_{+}(a_{1},...,a_{s}):=\sup\{\sum\limits_{n_{1},...,n_{s}}|T(x_{n_{1}}^{1},...,x_{n_{s}}^{s})|:x^{k}\in
\prod a_{k}(k\in \{1,...,s\})\}$ for every $a_{k}\in E_{k}^{+}\ (k\in\lbrace 1,...,s\rbrace)$. Like in the proof of Theorem 3.1 of \cite{BusvR}, one infers that $\bar{T}_{+}$ is additive and positively homogeneous in each variable separately. Therefore, by routine reasoning, $T_{+}$ uniquely extends to a positive $s_{\mathbb{R}}$-linear map $\bar{T}:\times_{k=1}^{s}E_{k\rho}\rightarrow F_{\rho}$, and subsequently to a positive $s_{\mathbb{C}}$-linear map $\bar{T}_{\mathbb{C}}:\times_{k=1}^{s}E_{k}\rightarrow F$. Then $\bar{T}_{\mathbb{C}}-T$ is also a
positive $s_{\mathbb{C}}$-linear map. By Theorem 4.13 there exists unique positive linear maps $\bar{T}_{\mathbb{C}}^{\bar{\otimes}}$ and $(\bar{T}_{\mathbb{C}}-T)^{\bar{\otimes}}$ from $\bar{\otimes}_{k=1}^{s}E_{k}$ into $F$ with $\bar{T}_{\mathbb{C}}=\bar{T}_{\mathbb{C}}^{\bar{\otimes}}\circ\bar{\otimes}$ and $(\bar{T}_{\mathbb{C}}-T)=(\bar{T}_{\mathbb{C}}-T)^{\bar{\otimes}}\circ\bar{\otimes}$. Then for $T^{\bar{\otimes}}:=\bar{T}_{\mathbb{C}}^{\bar{\otimes}}-(\bar{T}_{\mathbb{C}}-T)^{\bar{\otimes}}$ we have $T^{\bar{\otimes}}\circ\bar{\otimes}=T$.

(2) The map $\Phi:\mathcal{L}_{bv}(E_{1},...,E_{s};F)\rightarrow\mathcal{L}_{b}(\bar{\otimes}_{k=1}^{s}E_{k},F)$ is a $\mathbb{C}$-linear map. Let the map $T:\times_{k=1}^{s}E_{k}\rightarrow F$ be a positive $s$-linear map and let $\sum\limits_{j=1}^{n}f^{1}_{j}\otimes\dots\otimes f^{s}_{j}\in
\otimes_{k=1}^{s}E_{k}$ with $f^{k}_{j}\in E_{k}^{+}$ for every $j\in\lbrace 1,...,n\rbrace\ (k\in\lbrace
1,...,s\rbrace)$. Then $T^{\bar{\otimes}}(\sum\limits_{j=1}^{n}f^{1}_{j}\otimes\dots\otimes
f^{s}_{j})=\sum\limits_{j=1}^{n}T(f^{1}_{j},...,f^{s}_{j})\in F^{+}$. By relatively uniform density, $T^{\bar{\otimes}}$ is positive on $\bar{\otimes}_{k=1}^{s}E_{k}$. Therefore, $\Phi$ is positive.  By Theorem 4.13 we have that for every positive linear map $S:\bar{\otimes}_{k=1}^{s}E_{k}\rightarrow F$, there exists a unique positive map $T\in\mathcal{L}_{bv}(E_{1},...,E_{s};F)$ such that $T^{\bar{\otimes}}=S$. Therefore $\Phi$ and $\Phi^{-1}$ are ordered vector space isomorphisms with respect to the cones $\mathcal{L}^{+}_{bv}(E_{1},...,E_{s};F)$ and $\bigl(\mathcal{L}_{b}(\bar{\otimes}_{k=1}^{s}E_{k},F)\bigr)^{+}$. It follows from Lemma 5.6 that $\mathcal{L}_{bv}(E_{1},...,E_{s};F)$ is an Archimedean vector lattice over $\mathbb{C}$ with $\mathcal{L}^{+}_{bv}(E_{1},...,E_{s};F)$ as positive cone. Also by Lemma 5.6, $\Phi$ is a vector lattice isomorphism, and so $\mathcal{L}_{bv}(E_{1},...,E_{s};F)$ is Dedekind complete.

(3) Let $a_{k}\in E_{k}^{+}\ (k\in\lbrace 1,...,s\rbrace)$ and put $\theta\in[0,2\pi]$. Also let $T\in\mathcal{L}_{bv}(E_{1},...,E_{s};F)$, and let $\bar{T}$ be as in part (1) of this proof. Evidently, we have $\bar{T}(a_{1},...,a_{s})\geq\bigl(Re(e^{-i\theta}T)\bigr)(a_{1},...,a_{s})$ and thus $\bar{T}\geq\sup\lbrace Re(e^{-i\theta}T):\theta\in[0,2\pi]\rbrace=|T|$\ (see page 188 in \cite{Zan2} for this alternative definition of the modulus). On the other hand, if $x^{k}\in\prod a_{k}\ (k\in\lbrace 1,...,s\rbrace)$ then $\sum\limits_{n_{1},...,n_{s}}|T(x^{1}_{n_{1}},...,x^{s}_{n_{s}})|\leq\sum\limits_{n_{1},...,n_{s}}|T|(x^{1}_{n_{1}},...,x^{s}_{n_{s}})=|T|(a_{1},...,a_{s})$. Therefore $\bar{T}\leq|T|$.
\end{proof}


\begin{thebibliography}{99}

\bibitem{AB} C.D. Aliprantis, O. Burkinshaw, Positive Operators, Academic Press, Orlando, 1985.

\bibitem{Az} Y. Azouzi, Square mean closed real Riesz spaces, Ph.D. Dissert., Tunis, 2008.

\bibitem{AzBoBus} Y. Azzouzi, K. Boulabiar, G. Buskes, The de
Schipper formula and squares of Riesz spaces, Indag. Math. (N.S.) 17 (2006) no. 4, 479--496.

\bibitem{BS}  J. Bochnak, J. Siciak, Polynomials and multilinear
mappings in topological vector spaces, Studia Math. 39 (1971) 59--76.

\bibitem{BoBus}  K. Boulabiar, G. Buskes, Vector lattice powers:
f-algebras and functional calculus, Comm. Algebra 34 (2006) no. 4,
1435--1442.

\bibitem{BusSch} G. Buskes, C. Schwanke, Functional completions of Archimedean vector lattices, Unpublished results.

\bibitem{BusvR}  G. Buskes, A. van Rooij, Bounded variation and
tensor products of Banach lattices, Positivity 7 (2003) no. 1-2, 47--59.

\bibitem{BusvR2}  G. Buskes, A. van Rooij, Squares of Riesz spaces, Rocky Mountain J. Math. 31 (2001) no. 1, 45--56.

\bibitem{Fr}  D.H. Fremlin, Tensor products of Archimedean vector
lattices, Amer. J. Math. 94 (1972) 777--798.

\bibitem{Hag}  A.W. Hager, Some remarks on the tensor product of
function rings, Math. Z. 92 (1966) 210--224.

\bibitem{Kal}  N. J. Kalton, Hermitian operators on complex Banach
lattices and a problem of Garth Dales, J. Lond. Math. Soc. (2) 86 (2012) no. 3, 641--656.

\bibitem{Loane}  J. Loane, Polynomials on Riesz spaces, J. Math. Anal. Appl. 364 (2010) no. 1, 71--78.

\bibitem{Lotz}  H.P. Lotz, \"{U}ber das Spektrum positiver
Operatoren, Math. Z. 108 (1968) 15--32.

\bibitem{LuxZan1}  W.A.J. Luxemburg, A.C. Zaanen, Riesz Spaces, Vol. I., North-Holland, Amsterdam-London-New York, 1971.

\bibitem{LuxZan2}  W.A.J. Luxemburg, A.C. Zaanen, The linear
modulus of an order bounded linear transformation I, Nederl. Akad.
Wetensch. Proc. Ser. A 74 = Indag. Math. 33 (1971) 422--434.

\bibitem{MN}  P. Meyer-Nieberg, Banach Lattices, Springer-Verlag, Berlin-Heidelberg, 1991.

\bibitem{MW}  G. Mittelmeyer, M. Wolff, \"{U}ber den
Abolutbetrag auf komplexen Vektorverb\"{a}nden, Math. Z. 137 (1974) 87--92.

\bibitem{Pel}  A. Pe\l czy\'{n}ski, Some linear topological
properties of separable function algebras, Proc. Amer. Math. Soc. 18 (1967) 652--660.

\bibitem{Rieffel1}  M.A. Rieffel, A characterization of
commutative group algebras and measure algebras, Bull. Amer. Math. Soc. 69 (1963) 812--814.

\bibitem{Rieffel2}  M.A. Rieffel, A characterization of
commutative group algebras and measure algebras, Trans. Amer. Math. Soc.
116 (1965) 32--65.

\bibitem{Schae3}  H.H. Schaefer, Banach Lattices and Positive
Operators, Springer, Berlin-Heidelberg-New York, 1974.

\bibitem{Schae2}  H.H. Schaefer, Zur komplexen Erweiterung
linearer R\"{a}ume, Arch. Math. 10 (1959) 363--365.

\bibitem{Sh}  A.R. Schep, Factorization of positive multilinear
maps, Illinois J. Math. 28 (1984) no. 4, 579--591.

\bibitem{dS}  W.J.A. de Schipper, A note on the modulus of an order
bounded linear operator between complex vector lattices, Nederl. Akad.
Wetensch. Proc. Ser. A 76 = Indag. Math. 35
(1973) 355--367.

\bibitem{Vuz}  D. Vuza,  Sur les espaces vectoriels r\'{e}ticul\'{e}s
complexes, Rev. Roumaine Math. Pures Appl. 25 (1980) no. 4, 663--674.

\bibitem{Zan2}  A.C. Zaanen, Riesz spaces II, North-Holland,
Amsterdam, 1983.

\bibitem{vZ} G. van Zyl, Metrical aspects of the complexification of
tensor products and tensor norms, Ph.D. Dissert., Pretoria, 2009.

\end{thebibliography}
\end{document}